\documentclass[10pt,a4paper,titlepage]{article}
\usepackage[mathscr,mathcal]{eucal}
\usepackage{colortbl}
\usepackage{amsfonts}
\usepackage{amssymb}
\usepackage{amsmath}
\usepackage{xypic}
\newtheorem{theorem}{Theorem}[section]
\newtheorem{definition}[theorem]{Definition}
\newtheorem{ypoth}[theorem]{Assumptions}

\newtheorem{lemma}[theorem]{Lemma}

\newtheorem{corollary}[theorem]{Corollary}
\newtheorem{proposition}[theorem]{Proposition}

\newtheorem{cordef}[theorem]{Corollary-Definition}

\newcommand{\el}{\lambda}
\newcommand{\eb}{\beta}
\newtheorem{px}[theorem]{Example}

\newcommand{\ie}{\emph{i.e., }}

\begin{document}

\begin{center}
\textbf{\large Rouquier blocks of the cyclotomic Ariki-Koike algebras}\\
$ $\\
Maria Chlouveraki\\
\end{center}
$ $\\ \\
ABSTRACT. The definition of Rouquier for the families of characters of Weyl groups in terms of blocks of the associated Iwahori-Hecke algebra has made
possible the generalization of this notion to the complex reflection groups.
Here we give an algorithm for the determination of the ``Rouquier blocks''
of the cyclotomic Ariki-Koike algebras.

\section*{Introduction}

The work of G. Lusztig on the irreducible characters of reductive
groups over finite fields (cf. \cite{Lu1}) has displayed the important role of the
``families of characters'' of the Weyl groups concerned. 
More recent results of Gyoja $\cite{Gy}$ and Rouquier
$\cite{Rou}$ have made possible the definition of a substitute for
families of characters which can be applied to all complex
reflection groups. In particular, Rouquier has shown  that the families
of characters of a Weyl group $W$ are exactly the blocks of the Iwahori-Hecke algebra of $W$ over a suitable
coefficient ring, the ``Rouquier ring''. This definition generalizes without problem to the
cyclotomic Hecke algebras of all complex reflection groups. Ever since, we have been interested
in the determination of the ``Rouquier blocks'' of the cyclotomic Hecke algebras of complex reflection groups.

In \cite{BK}, Brou\'{e} and Kim presented an algorithm for the determination of the Rouquier blocks for the cyclotomic Hecke algebras of the groups $G(d,1,r)$ and $G(d,d,r)$. Later Kim \cite{Kim} used the same algorithm to determine the Rouquier blocks for the group $G(de,e,r)$. Moreover, the ``Rouquier
blocks'' of the ``spetsial'' cyclotomic Hecke algebra of many
exceptional complex reflection groups have been determined by Malle
and Rouquier in \cite{MaRo}. Finally, in \cite{Chlou}, we've determined the Rouquier blocks of the cyclotomic Hecke algebras of all exceptional complex reflection groups. 

However,  it was recently realized that the algorithm given in \cite{BK} works only in the case where 
$d$ is the power of a prime number. The aim of this paper is to give a complete description of the Rouquier blocks of  the cyclotomic Ariki-Koike algebras of the group $G(d,1,r)$. In order to achieve that, we use the theory of ``essential hyperplanes'' introduced in \cite{Chlou}. According to this theory, the Rouquier blocks of the cyclotomic Hecke algebras of any complex reflection group depend on numerical data determined by the generic Hecke algebra, the ``essential hyperplanes'' of the group. Thanks to Theorem $\ref{main theorem}$, it suffices to study the blocks of the generic Hecke algebra in a finite number of cases in order to obtain the Rouquier blocks for all cyclotomic Hecke algebras.

An algorithm for the blocks of the Ariki-Koike algebras of $G(d,1,r)$ over any field has been given by Lyle and Mathas in \cite{LyMa}. This algorithm can be applied to give us the Rouquier blocks of the cyclotomic Ariki-Koike algebras and we use it to obtain a characterization in the combinatorial terms used in \cite{BK}. Our main result is Theorem $\ref{RB}$, which determines completely the Rouquier blocks of the cyclotomic Ariki-Koike algebras. The most important consequence is that we can obtain the Rouquier blocks of a cyclotomic Ariki-Koike algebra of $G(d,1,r)$ from the families of characters of the Weyl groups of type $B_n$, $n \leq r$, already determined by Lusztig. This result can also be deduced from the Morita equivalences established by Dipper and Mathas in \cite{DiMa}. 
Moreover, we show that the Rouquier blocks in the important case of the ``spetsial'' cyclotomic Hecke algebra are the ones given by the algorithm of \cite{BK}.

Finally, in the case of the Weyl groups, Lusztig attaches to every irreducible character two integers, denoted by $a$ and $A$, and shows (cf. \cite{Lu2}, $3.3$ and $3.4$) that they are constant on the families. 
In an analogue way, we can define integers $a$ and $A$ attached to every irreducible character of a cyclotomic Hecke algebra of a complex reflection group. Proposition $\ref{a}$ completes the proof of the result by Brou\'{e}-Kim (cf. \cite{BK}, $3.18$) that the integers $a$ and $A$ are constant on the Rouquier blocks of $G(d,1,r)$. The same result has been obtained by the author for  the exceptional complex reflection groups in \cite{ChDeg}.
\\
\\
\textbf{Acknowledgements}. I would like to thank C\'{e}dric Bonnaf\'{e} and Nicolas Jacon for our fruitful conversations and their useful remarks.  I would also like to thank the Mathematical Sciences Research Institute of Berkeley and the Ecole Polytechnique F\'{e}d\'{e}rale de Lausanne for their financial support.

\section{Blocks and Symmetric algebras}

For the proofs of the results not given in this section, the reader may refer to \cite{BK} or the second chapter of \cite{Chlou}. 

\subsection{Generalities}

Let us assume that $\mathcal{O}$ is a commutative integral domain with field of fractions
  $F$ and $A$ is an $\mathcal{O}$-algebra, free and finitely generated as an
  $\mathcal{O}$-module. We denote by $ZA$ the center of $A$.
 
\begin{definition}\label{blocks}
The block-idempotents (blocks) of $A$ are the primitive idempotents
of $ZA$.
\end{definition}

Let $K$ be a field extension of $F$ and suppose that the $K$-algebra $KA:=K \otimes_\mathcal{O}A$ is
semisimple. Then $KA$ is isomorphic, by assumption, to a direct
product of simple algebras:
$$KA \simeq \prod_{\chi \in \mathrm{Irr}(KA)} M_\chi,$$
where $\mathrm{Irr}(KA)$ denotes the set of irreducible characters
of $KA$ and $M_\chi$ is a simple $K$-algebra.

For all $\chi \in \mathrm{Irr}(KA)$, we denote by $\pi_\chi:KA
\twoheadrightarrow M_\chi$ the projection onto the $\chi$-factor and
by $e_\chi$ the element of $KA$ such that $$\pi_{\chi'}(e_\chi)=
  \left\{
  \begin{array}{ll}
    1_{M_\chi}, & \hbox{if $\chi=\chi'$,} \\
    0, & \hbox{if $\chi \neq \chi'$.}
  \end{array}
\right.$$

We have the following result which provides a relation between the blocks of the algebra $KA$ and the blocks of the algebra $A$.

\begin{theorem}\label{minimality of blocks}\
\begin{enumerate}
  \item We have $1=\sum_{\chi \in \mathrm{Irr}(KA)}e_\chi$
   and the set $\{e_\chi\}_{\chi \in \mathrm{Irr}(KA)}$ is the set of all the blocks of the algebra $KA$.
  \item There exists a unique partition $\mathrm{Bl}(A)$ of
  $\mathrm{Irr}(KA)$ such that
  \begin{description}
    \item[(a)] For all $B \in \mathrm{Bl}(A)$, the idempotent
    $e_B:=\sum_{\chi \in B}e_\chi$ is a block of $A$.
    \item[(b)] We have $1=\sum_{B \in \mathrm{Bl}(A)}e_B$ and for
    every central idempotent $e$ of $A$, there exists a subset
    $\mathrm{Bl}(A,e)$ of $\mathrm{Bl}(A)$ such that
    $$e=\sum_{B \in \mathrm{Bl}(A,e)}e_B.$$
    In particular the set $\{e_B\}_{B \in \mathrm{Bl}(A)}$ is the set of all the blocks of $A$.
  \end{description}
\end{enumerate}
\end{theorem}

\begin{remarks}\
\emph{\begin{itemize}
  \item If $\chi \in B$ for some $B \in \mathrm{Bl}(A)$, we say that
  ``$\chi$ belongs to the block $e_B$''.
  \item For all $B \in \mathrm{Bl}(A)$, we have
  $$KAe_B \simeq \prod_{\chi \in B}M_\chi.$$
\end{itemize}}
\end{remarks}

From now on, we make the following assumptions

\begin{ypoth}\label{properties of the ring}\
\begin{description}
  \item[(int)] The ring $\mathcal{O}$ is a Noetherian and integrally
  closed domain with field of fractions $F$ and $A$ is an
  $\mathcal{O}$-algebra which is free and finitely generated as an
  $\mathcal{O}$-module.
  \item[(spl)] The field $K$ is a finite Galois extension of $F$ and
  the algebra $KA$ is split (i.e., for every simple $KA$-module $V$, $\mathrm{End}_{KA}(V) \simeq K$) semisimple.
\end{description}
\end{ypoth}

We denote by $\mathcal{O}_K$ the integral closure of $\mathcal{O}$
in $K$.

\subsection{Blocks and integral closure}

The Galois group $\mathrm{Gal}(K/F)$ acts on $KA=K
\otimes_{\mathcal{O}} A$ (viewed as an $F$-algebra) as follows: if
$\sigma \in \mathrm{Gal}(K/F)$ and $\lambda \otimes a \in KA$, then
$\sigma(\lambda \otimes a):=\sigma(\lambda) \otimes a$.

If $V$ is a $K$-vector space and $\sigma \in \mathrm{Gal}(K/F)$, we
denote by $^\sigma V$ the $K$-vector space defined on the additive
group $V$ with multiplication $\lambda.v:=\sigma^{-1}(\lambda)v$ for
all $\lambda \in K$ and $v \in V$. If $\rho:KA \rightarrow
\mathrm{End}_K(V)$ is a representation of the $K$-algebra $KA$, then
its composition with the action of $\sigma^{-1}$ is also a
representation $^\sigma \rho: KA \rightarrow \mathrm{End}_K(^\sigma
V)$:
$$\diagram KA \rto^{\sigma^{-1}} &KA \rto^{\rho}
&\mathrm{End}_K(V). \enddiagram$$

We denote by $^\sigma \chi$ the character of $^\sigma \rho$ and we
define the action of $\mathrm{Gal}(K/F)$ on $\mathrm{Irr}(KA)$ as
follows: if $\sigma \in \mathrm{Gal}(K/F)$ and $\chi \in
\mathrm{Irr}(KA)$, then $$\sigma(\chi):={}^\sigma\!\chi = \sigma
\circ \chi \circ \sigma^{-1}.$$ This operation induces an action of
$\mathrm{Gal}(K/F)$ on the set of blocks of $KA$:
$$\sigma(e_\chi)=e_{^\sigma \chi} \textrm{ for all } \sigma \in
\mathrm{Gal}(K/F), \chi \in \mathrm{Irr}(KA).$$

Hence, the group $\mathrm{Gal}(K/F)$ acts on the set of idempotents
of $Z\mathcal{O}_KA$ and thus on the set of blocks of
$\mathcal{O}_KA$. Since $F \cap \mathcal{O}_K = \mathcal{O}$, the
idempotents of $ZA$ are the idempotents of $Z\mathcal{O}_KA$ which
are fixed by the action of $\mathrm{Gal}(K/F)$. As a consequence,
the primitive idempotents of $ZA$ are sums of the elements of the
orbits of $\mathrm{Gal}(K/F)$ on the set of primitive idempotents of
$Z\mathcal{O}_KA$. Thus, the blocks of $A$ are in bijection with the
orbits of $\mathrm{Gal}(K/F)$ on the set of blocks of
$\mathcal{O}_KA$. The following proposition is just a reformulation
of this result.

\begin{proposition}\label{Galois action on integral closure}\
\begin{enumerate}
  \item Let $B$ be a block of $A$ and $B'$ a block of
  $\mathcal{O}_KA$ contained in $B$. If $\mathrm{Gal}(K/F)_{B'}$
  denotes the stabilizer of $B'$ in $\mathrm{Gal}(K/F)$, then
  $$B=\bigcup_{\sigma \in
  \mathrm{Gal}(K/F)/\mathrm{Gal}(K/F)_{B'}}\sigma(B')
  \,\textrm{ i.e., }\,
  e_B=\sum_{\sigma \in
  \mathrm{Gal}(K/F)/\mathrm{Gal}(K/F)_{B'}}\sigma(e_{B'}).$$
  \item Two characters $\chi,\psi \in \mathrm{Irr}(KA)$ are in
  the same block of $A$ if and only if there exists $\sigma \in \mathrm{Gal}(K/F)$
  such that $\sigma(\chi)$ and $\psi$ belong to the same block of
  $\mathcal{O}_KA$.
\end{enumerate}
\end{proposition}\
\begin{remark}\emph{ For all $\chi \in B'$, we have
$\mathrm{Gal}(K/F)_\chi \subseteq \mathrm{Gal}(K/F)_{B'}.$}
\end{remark}\
\\

The assertion (2) of the proposition above allows us to transfer the
problem of the classification of the blocks of $A$ to that of the
classification of the blocks of $\mathcal{O}_KA$.

\subsection{Blocks and prime ideals}

We denote by $\mathrm{Spec}(\mathcal{O})$ the set of prime ideals
of $\mathcal{O}$. Since $\mathcal{O}$ is Noetherian and
integrally closed, we have
$$\mathcal{O}=\bigcap_{\mathfrak{p} \in \mathrm{Spec}(\mathcal{O})}
\mathcal{O}_\mathfrak{p},$$ where $\mathcal{O}_\mathfrak{p}:=\{x \in
F\,|\,(\exists a \in \mathcal{O}-\mathfrak{p})(ax \in
\mathcal{O})\}$ is the localization of $\mathcal{O}$ at
$\mathfrak{p}$. 

Let $\mathfrak{p}$ be a prime ideal of $\mathcal{O}$ and
$\mathcal{O}_\mathfrak{p}A:=\mathcal{O}_\mathfrak{p}
\otimes_{\mathcal{O}}A$. The blocks of $\mathcal{O}_\mathfrak{p}A$
are the ``$\mathfrak{p}$-blocks of $A$''. If $\chi,\psi \in
\mathrm{Irr}(KA)$ belong to the same block of
$\mathcal{O}_\mathfrak{p}A$, we write $\chi \sim_\mathfrak{p} \psi$.

\begin{proposition}\label{p-blocks}
Two characters $\chi,\psi \in \mathrm{Irr}(KA)$ belong to the same
block of $A$ if and only if there exist a finite sequence
$\chi_0,\chi_1,\ldots,\chi_n \in \mathrm{Irr}(KA)$ and a finite
sequence $\mathfrak{p}_1,\ldots,\mathfrak{p}_n \in
\mathrm{Spec}(\mathcal{O})$ such that
\begin{itemize}
  \item $\chi_0=\chi$ and $\chi_n=\psi$,
  \item for all $j$ $(1\leq j \leq n)$, $\chi_{j-1}
  \sim_{\mathfrak{p}_j} \chi_j$.
\end{itemize}
\end{proposition}

\subsection{Blocks and residue blocks}

Let $\mathfrak{p}$ be a maximal ideal of $\mathcal{O}$ and set $k_{\mathfrak{p}}:=\mathcal{O}/\mathfrak{p}$ its residue field. If $\mathcal{O}_\mathfrak{p}$ is the localization of $\mathcal{O}$ at $\mathfrak{p}$, then $k_{\mathfrak{p}}$ is also the residue field of $\mathcal{O}_\mathfrak{p}$. The natural surjection $\pi_\mathfrak{p}:\mathcal{O}_\mathfrak{p} \twoheadrightarrow k_\mathfrak{p}$ extends to a morphism 
$\pi_\mathfrak{p}:\mathcal{O}_\mathfrak{p}A \twoheadrightarrow k_\mathfrak{p}A$, which in turn induces a morphism $$\pi_\mathfrak{p}:Z\mathcal{O}_\mathfrak{p}A \rightarrow Zk_\mathfrak{p}A. $$ 
The following lemma will serve for the proof of proposition $\ref{residue field}$.

\begin{lemma}\label{central lifts to central}
Let  $e$ be an idempotent of $\mathcal{O}_\mathfrak{p}A$ whose image $\bar{e}$ in
$k_\mathfrak{p}A$ is central. Then $e$ is
central.
\end{lemma}
\begin{apod}{Set $R:=\mathcal{O}_\mathfrak{p}A$. Since $\bar{e}$ is central, we have
$\bar{e}k_\mathfrak{p}A(1-\bar{e})=(1-\bar{e})k_\mathfrak{p}A\bar{e}=\{0\},$ \ie
$eR(1-e) \subseteq \mathfrak{p}R \,\,\textrm{
and }\,\, (1-e)Re \subseteq \mathfrak{p}R.$
Since $e$ and $(1-e)$ are idempotents, we get
$eR(1-e) \subseteq \mathfrak{p}eR(1-e) \,\,\textrm{ and
}\,\, (1-e)Re \subseteq
\mathfrak{p}(1-e)Re.$  By Nakayama's lemma,
$eR(1-e)=(1-e)Re=\{0\}.$
Thus, from
$R=eRe \oplus eR(1-e) \oplus (1-e)Re \oplus (1-e)R(1-e)$
we deduce that
$R=eRe \oplus (1-e)R(1-e)$
and consequently, $e$ is central.}
\end{apod}
\begin{proposition}\label{residue field}
If $\mathcal{O}_\mathfrak{p}$ is a discrete valuation ring and $K=F$, then the morphism 
$$\pi_\mathfrak{p}:Z\mathcal{O}_\mathfrak{p}A \rightarrow Zk_\mathfrak{p}A$$ 
induces a bijection between the set of blocks of $\mathcal{O}_\mathfrak{p}A$ and the set of blocks of
$k_\mathfrak{p}A$.
\end{proposition}
\begin{apod}{From now on, the symbol $ $ $\widehat{}$ $ $ will stand for
$\mathfrak{p}$-adic completion. It is immediate that $\pi_\mathfrak{p}$ sends a block of $\mathcal{O}_\mathfrak{p}A$  to a sum of blocks of $k_\mathfrak{p}A$.
Now let $\bar{e}$ be a block of $k_\mathfrak{p}A$.
By the theorems of lifting idempotents (see \cite{The}, Thm.3.2) and
the lemma above, $\bar{e}$ is lifted to a sum of central
primitive idempotents in $\hat{\mathcal{O}}_{\mathfrak{p}}A$.
However, by the fact that $KA$ is
split semisimple, we have that the blocks of
$\hat{\mathcal{O}}_{\mathfrak{p}}A$ belong to
$KA$. But $K \cap
\hat{\mathcal{O}}_{\mathfrak{p}}=\mathcal{O}_{\mathfrak{p}}$ (cf., for example, \cite{Na}, 18.4)
and $\mathcal{O}_{\mathfrak{p}}A \cap
Z \hat{\mathcal{O}}_{\mathfrak{p}}A \subseteq Z \mathcal{O}_{\mathfrak{p}}A$. Therefore, $\bar{e}$ is lifted to
a sum of blocks in $\mathcal{O}_{\mathfrak{p}}A$ and this provides
the block bijection.
\end{apod}}

\subsection{Symmetric algebras}

Let $\mathcal{O}$ be a ring and let $A$ be an $\mathcal{O}$-algebra.
Suppose again that the assumptions $\ref{properties of the ring}$
are satisfied.

\begin{definition}\label{trace function}
A trace function on $A$ is an $\mathcal{O}$-linear map $t:A
\rightarrow \mathcal{O}$ such that $t(ab)=t(ba)$ for all $a,b \in
A$.
\end{definition}

\begin{definition}\label{symmetric algebra}
We say that a trace function $t:A \rightarrow \mathcal{O}$ is a
symmetrizing form on $A$ or that $A$ is a symmetric algebra if the
morphism
$$\hat{t}:A \rightarrow \mathrm{Hom}_\mathcal{O}(A,\mathcal{O}),\,\,
  a \mapsto (x \mapsto \hat{t}(a)(x):=t(ax))$$
is an isomorphism of $A$-modules-$A$.
\end{definition}

\begin{px}\label{symmetrizing form of the group algebra}
\small{\emph{In the case where $\mathcal{O}=\mathbb{Z}$ and
$A=\mathbb{Z}[G]$
 ($G$ a finite group), we can define the following symmetrizing form
 (``canonical'')
 on $A$
$$t:\mathbb{Z}[G] \rightarrow \mathbb{Z}, \,\,\, \sum_{g \in G}a_g g \mapsto a_1,$$
where $a_g \in \mathbb{Z}$ for all $g \in G$.}}
\end{px}

If $\tau:A \rightarrow \mathcal{O}$ is a linear form, we denote by
$\tau^\vee$ its inverse image by the isomorphism $\hat{t}$, \ie
$\tau^\vee$ is the element of $A$ such that
$$t(\tau^\vee a)=\tau(a) \textrm{ for all } a \in A.$$

The element $\tau^\vee$ has the following properties (cf., for example,  \cite{GePf}, \S 7.1):

\begin{lemma}\label{tau^vee}\
\begin{enumerate}
 \item $\tau$ is a trace function if and only if $\tau^\vee \in ZA$.
 \item Let $(e_i)_{i \in I}$ be a basis of $A$ over
$\mathcal{O}$ and $(e_i')_{i \in I}$ is the dual basis with respect
to $t$ (i.e., $t(e_ie_j')=\delta_{ij}$).
 We have $\tau^\vee=\sum_i \tau(e_i')e_i=\sum_i \tau(e_i)e_i'$ and
 more generally, for all $a \in A$, we have
 $\tau^\vee a=\sum_i \tau(e_i'a)e_i=\sum_i \tau(e_ia)e_i'$.
\end{enumerate}
\end{lemma}

\subsection{Schur elements}

If $A$ is a symmetric algebra with a symmetrizing form $t$, we
obtain a symmetrizing form $t^K$ on $KA$ by extension of scalars.
Every irreducible character $\chi \in \mathrm{Irr}(KA)$ is a trace
function on $KA$ and thus we can define $\chi^\vee \in ZKA$.
Since $KA$ is a split semisimple $K$-algebra, we have that
$KA \simeq \prod_{\chi \in \mathrm{Irr}(KA)} M_\chi,$
where $M_\chi$ is a matrix algebra isomorphic to
$\mathrm{Mat}_{\chi(1)}(K)$. The map $\pi_\chi: KA
\twoheadrightarrow M_\chi$, restricted to $ZKA$, defines a map
$\omega_\chi:ZKA \twoheadrightarrow K$.

\begin{definition}\label{Schur element}
For all $\chi \in \mathrm{Irr}(KA)$, we call Schur element of $\chi$
with respect to $t$ and denote by $s_\chi$ the element of $K$
defined by $$s_\chi:=\omega_\chi(\chi^\vee).$$
\end{definition}

The following property of the Schur elements is proven in \cite{GePf}, \S 7.2.

\begin{proposition}\label{Schur element belongs to the integral closure}
For all $\chi \in \mathrm{Irr}(KA)$, $s_\chi \in {\mathcal{O}_K}^*$.
\end{proposition}

\begin{px}\label{Schur elements of the group algebra}
\small{\emph{Let $\mathcal{O}:=\mathbb{Z}$, $A:=\mathbb{Z}[G]$
 ($G$ a finite group) and $t$ the canonical symmetrizing form. If $K$ is an algebraically closed field of
 characteristic 0, then $KA$ is a split semisimple algebra and
 $s_\chi=|G|/\chi(1)$ for all $\chi \in \mathrm{Irr}(KA)$. Because
 of the integrality of the Schur elements, we must have
 $|G|/\chi(1) \in \mathbb{Z}=\mathbb{Z}_K \cap \mathbb{Q}$ for all $\chi \in
 \mathrm{Irr}(KA)$. Thus, we have shown that $\chi(1)$ divides $|G|$.}}
\end{px}

The following properties of the Schur elements can be derived easily
from the above (see also
\cite{Bro},\cite{Ge},\cite{GePf},\cite{GeRo},\cite{BMM2}).

\begin{proposition}\label{schur elements and idempotents}\
\begin{enumerate}
  \item We have
  $$t=\sum_{\chi \in \mathrm{Irr}(KA)}\frac{1}{s_\chi}\chi.$$
  \item For all $\chi \in \mathrm{Irr}(KA)$, the central primitive
  idempotent associated with $\chi$ is
  $$e_\chi=\frac{1}{s_\chi}\chi^\vee.$$
 \end{enumerate}
\end{proposition}

\section{Hecke algebras of complex reflection groups}

\subsection{Generic Hecke algebras}

Let $\mu_\infty$ be the group of all the roots of unity in
$\mathbb{C}$ and $K$ a number field contained in
$\mathbb{Q}(\mu_\infty)$. We denote by $\mu(K)$ the group of all the
roots of unity of $K$. For every integer $d>1$, we set
$\zeta_d:=\mathrm{exp}(2\pi i/d)$ and denote by $\mu_d$ the group of
all the $d$-th roots of unity. 

Let $V$ be a $K$-vector space of
finite dimension $r$. Let $W$ be a finite subgroup of $\mathrm{GL}(V)$ generated by
(pseudo-)reflections acting irreducibly on $V$. Let us denote by $\mathcal{A}$ the set of the
reflecting hyperplanes of $W$. We set $\mathcal{M} := \mathbb{C} \otimes V - 
\bigcup_{H \in \mathcal{A}} \mathbb{C} \otimes H$. For $x_0 \in \mathcal{M}$, 
let $P:=\Pi_1(\mathcal{M},x_0)$ and $B:=\Pi_1(\mathcal{M}/W,x_0)$. Then there exists a short exact 
sequence (cf. \cite{BMR}, \S 2B):
$$\{1\}\rightarrow P \rightarrow B \rightarrow W \rightarrow\{1\}.$$
We denote by $\tau$ the central element of $P$ defined by the loop
$$[0,1] \rightarrow \mathcal{M}, \,\,\,\,\,\,\, t \mapsto \mathrm{exp}(2\pi it)x_0.$$ 

For every orbit $\mathcal{C}$ of $W$ on $\mathcal{A}$, we denote by
$e_{\mathcal{C}}$ the common order of the subgroups $W_H$, where $H$
is any element of $\mathcal{C}$ and $W_H$ the subgroup formed by $\mathrm{id}_V$
and all the reflections fixing the hyperplane $H$.

We choose a set of indeterminates
$\textbf{u}=(u_{\mathcal{C},j})_{(\mathcal{C} \in
\mathcal{A}/W)(0\leq j \leq e_{\mathcal{C}}-1)}$ and we denote by
$\mathbb{Z}[\textbf{u},\textbf{u}^{-1}]$ the Laurent polynomial ring
in all the indeterminates $\textbf{u}$. We define the \emph{generic
Hecke algebra} $\mathcal{H}$ of $W$ to be the quotient of the group
algebra $\mathbb{Z}[\textbf{u},\textbf{u}^{-1}]B$ by the ideal
generated by the elements of the form
$$(\textbf{s}-u_{\mathcal{C},0})(\textbf{s}-u_{\mathcal{C},1}) \ldots (\textbf{s}-u_{\mathcal{C},e_{\mathcal{C}}-1}),$$
where $\mathcal{C}$ runs over the set $\mathcal{A}/W$ and
$\textbf{s}$ runs over the set of monodromy generators around the
images in $\mathcal{M}/W$ of the elements of the hyperplane
orbit $\mathcal{C}$.

\begin{px}
\small{\emph{Let $W:=G_4=<s,t \,|\, \,sts=tst, s^3=t^3=1>$. Then $s$
and $t$ are conjugate in $W$ and their reflecting hyperplanes belong
to the same orbit in $\mathcal{A}/W$. The generic Hecke algebra of
$W$ can be presented as follows
$$\begin{array}{rll}
   \mathcal{H}(G_4)=<S,T \,\,|&STS=TST, &(S-u_0)(S-u_1)(S-u_2)=0, \\
                                &         &(T-u_0)(T-u_1)(T-u_2)=0>.
  \end{array}$$}}
\end{px}

We make some assumptions for the algebra $\mathcal{H}$. Note that
they have been verified for all but a finite number of irreducible
complex reflection groups (\cite{BMM2}, remarks before 1.17, $\S$ 2;
\cite{GIM}).

\begin{ypoth}\label{ypo}
The algebra $\mathcal{H}$ is a free
$\mathbb{Z}[\textbf{\emph{u}},\textbf{\emph{u}}^{-1}]$-module of
rank $|W|$. Moreover, there exists a linear form
$t:\mathcal{H}\rightarrow
\mathbb{Z}[\textbf{\emph{u}},\textbf{\emph{u}}^{-1}]$ with the
following properties:
\begin{enumerate}
    \item $t$ is a symmetrizing form for $\mathcal{H}$.
    \item Via the specialization $u_{\mathcal{C},j} \mapsto
     \zeta_{e_\mathcal{C}}^j$, the form $t$ becomes the canonical
     symmetrizing form on the group algebra $\mathbb{Z}W$.
    \item If we denote by $\alpha \mapsto \alpha^*$ the automorphism of
     $\mathbb{Z}[\emph{\textbf{u}},\emph{\textbf{u}}^{-1}]$ consisting of the
     simultaneous inversion of the indeterminates, then for all $b \in B$, we
     have
          $$t(b^{-1})^*=\frac{t(b\tau)}{t(\tau)}.$$
\end{enumerate}
\end{ypoth}

We know that the form $t$ is unique (\cite{BMM2}, 2.1). From now on,
let us suppose that the assumptions $\ref{ypo}$ are satisfied. Then
we have the following result by G.Malle (\cite{Ma4}, 5.2).

\begin{theorem}\label{Semisimplicity Malle}
Let $\textbf{\emph{v}}=(v_{\mathcal{C},j})_{(\mathcal{C} \in
\mathcal{A}/W)(0\leq j \leq e_{\mathcal{C}}-1)}$ be a set of
$\sum_{\mathcal{C} \in \mathcal{A}/W}e_{\mathcal{C}}$ indeterminates
such that, for every $\mathcal{C},j$, we have
$v_{\mathcal{C},j}^{|\mu(K)|}=\zeta_{e_\mathcal{C}}^{-j}u_{\mathcal{C},j}$.
Then the $K(\textbf{\emph{v}})$-algebra
$K(\textbf{\emph{v}})\mathcal{H}$ is split semisimple.
\end{theorem}

By ``Tits' deformation theorem'' (cf., for example, \cite{BMM2}, 7.2), it follows
that the specialization $v_{\mathcal{C},j}\mapsto 1$ induces a
bijection $\chi \mapsto \chi_{\textbf{v}}$ from the
set $\mathrm{Irr}(K(\textbf{v})\mathcal{H})$ of absolutely
irreducible characters of $K(\textbf{v})\mathcal{H}$ to the set
$\mathrm{Irr}(W)$ of absolutely irreducible characters of $W$, such that the
following diagram is commutative $$\begin{array}{rccc}
  \chi_\textbf{v} : & \mathcal{H} & \rightarrow & \mathbb{Z}_K[\textbf{v},\textbf{v}^{-1}] \\
  & \downarrow &  & \downarrow \\
  \chi: & \mathbb{Z}_KW &\rightarrow &\mathbb{Z}_K .
\end{array}$$
$ $\\

The following result concerning the form of the Schur elements associated
with the irreducible characters of $K(\textbf{v})\mathcal{H}$ is proved in \cite{Chlou}, Thm.3.2.5, using case by case analysis.

\begin{theorem}\label{Schur element generic}
The Schur element $s_\chi(\textbf{\emph{v}})$ associated with the
character $\chi_{\textbf{\emph{v}}}$ of
$K(\textbf{\emph{v}})\mathcal{H}$ is an element of
$\mathbb{Z}_K[\textbf{\emph{v}},\textbf{\emph{v}}^{-1}]$ of the form
$$s_\chi({\textbf{\emph{v}}})=\xi_\chi N_\chi \prod_{i \in I_\chi} \Psi_{\chi,i}(M_{\chi,i})^{n_{\chi,i}}$$
where
\begin{itemize}
    \item $\xi_\chi$ is an element of $\mathbb{Z}_K$,
    \item $N_\chi= \prod_{\mathcal{C},j} v_{\mathcal{C},j}^{b_{\mathcal{C},j}}$ is a monomial in $\mathbb{Z}_K[\textbf{\emph{v}},\textbf{\emph{v}}^{-1}]$
          such that $\sum_{j=0}^{e_\mathcal{C}-1}b_{\mathcal{C},j}=0$
          for all $\mathcal{C} \in \mathcal{A}/W$,
    \item $I_\chi$ is an index set,
    \item $(\Psi_{\chi,i})_{i \in I_\chi}$ is a family of $K$-cyclotomic polynomials in one variable
           (i.e., minimal polynomials of the roots of unity over $K$),
    \item $(M_{\chi,i})_{i \in I_\chi}$ is a family of monomials in $\mathbb{Z}_K[\textbf{\emph{v}},\textbf{\emph{v}}^{-1}]$
          and if $M_{\chi,i} = \prod_{\mathcal{C},j} v_{\mathcal{C},j}^{a_{\mathcal{C},j}}$,
          then $\textrm{\emph{gcd}}(a_{\mathcal{C},j})=1$
          and $\sum_{j=0}^{e_\mathcal{C}-1}a_{\mathcal{C},j}=0$
          for all $\mathcal{C} \in \mathcal{A}/W$,
    \item ($n_{\chi,i})_{i \in I_\chi}$ is a family of positive integers.
\end{itemize}
This factorization is unique in $K[\textbf{\emph{v}},\textbf{\emph{v}}^{-1}]$. Moreover, the monomials
 $(M_{\chi,i})_{i \in I_\chi}$ are unique up to inversion, whereas the coefficient $\xi_\chi$ is unique up to multiplication by a root of unity.
\end{theorem}
\begin{remark}
\emph{The bijection $\mathrm{Irr}(K(\textbf{v})\mathcal{H}) \leftrightarrow \mathrm{Irr}(W)$, $\chi_{\textbf{v}} \mapsto
\chi$  implies that the
specialization $v_{\mathcal{C},j}\mapsto 1$ sends
$s_{\chi_{\textbf{v}}}$ to $|W|/\chi(1)$ (which is the Schur element
of $\chi$ in the group algebra with respect to the canonical
symmetrizing form). }
\end{remark}\\

Let $A:=\mathbb{Z}_K[\textbf{v},\textbf{v}^{-1}]$ and $\mathfrak{p}$ be a prime ideal of $\mathbb{Z}_K$. 

\begin{definition}\label{p-essential monomial}
Let  $M = \prod_{\mathcal{C},j} v_{\mathcal{C},j}^{a_{\mathcal{C},j}}$ be a monomial in $A$
such that $\textrm{\emph{gcd}}(a_{\mathcal{C},j})=1$. We say that $M$ is $\mathfrak{p}$-essential
for a character $\chi \in \mathrm{Irr}(W)$, if there exists a $K$-cyclotomic polynomial $\Psi$ such that
\begin{itemize}
\item $\Psi(M)$ divides $s_\chi(\textbf{\emph{v}})$.
\item $\Psi(1)  \in  \mathfrak{p}$.
\end{itemize}
We say that $M$ is $\mathfrak{p}$-essential 
for $W$, if there exists a character $\chi \in \mathrm{Irr}(W)$ such that
$M$ is $\mathfrak{p}$-essential for $\chi$.
\end{definition}

The following proposition (\cite{Chlou}, Prop.3.2.6) gives a characterization of $\mathfrak{p}$-essential monomials, which plays an essential role in the proof of theorem $\ref{main theorem}$. 

\begin{proposition}\label{p-essential}
Let  $M = \prod_{\mathcal{C},j} v_{\mathcal{C},j}^{a_{\mathcal{C},j}}$ be a monomial in $A$
such that $\textrm{\emph{gcd}}(a_{\mathcal{C},j})=1$. We set  $\mathfrak{q}_M:=(M-1)A +\mathfrak{p}A$.
Then
\begin{enumerate}
\item The ideal $\mathfrak{q}_M$ is a prime ideal of $A$. 
\item $M$ is $\mathfrak{p}$-essential for $\chi \in \mathrm{Irr}(W)$ if and only if
$s_\chi(\textbf{\emph{v}})/\xi_\chi \in \mathfrak{q}_M$.
\end{enumerate}
\end{proposition}

\subsection{Cyclotomic Hecke algebras}

Let $y$ be an indeterminate. We set $x:=y^{|\mu(K)|}.$

\begin{definition}\label{specialization}
A cyclotomic specialization of $\mathcal{H}$ is a
$\mathbb{Z}_K$-algebra morphism $\phi:
\mathbb{Z}_K[\textbf{\emph{v}},\textbf{\emph{v}}^{-1}]\rightarrow
\mathbb{Z}_K[y,y^{-1}]$ with the following properties:
\begin{itemize}
  \item $\phi: v_{\mathcal{C},j} \mapsto y^{n_{\mathcal{C},j}}$ where
  $n_{\mathcal{C},j} \in \mathbb{Z}$ for all $\mathcal{C}$ and $j$.
  \item For all $\mathcal{C} \in \mathcal{A}/W$, if $z$ is another
  indeterminate, the element of $\mathbb{Z}_K[y,y^{-1},z]$ defined by
  $$\Gamma_\mathcal{C}(y,z):=\prod_{j=0}^{e_\mathcal{C}-1}(z-\zeta_{e_\mathcal{C}}^jy^{n_{\mathcal{C},j}})$$
  is invariant by the action of $\textrm{\emph{Gal}}(K(y)/K(x))$.
\end{itemize}
\end{definition}

If $\phi$ is a cyclotomic specialization of $\mathcal{H}$,
the corresponding \emph{cyclotomic Hecke algebra} is the
$\mathbb{Z}_K[y,y^{-1}]$-algebra, denoted by $\mathcal{H}_\phi$,
which is obtained as the specialization of the
$\mathbb{Z}_K[\textbf{v},\textbf{v}^{-1}]$-algebra $\mathcal{H}$ via
the morphism $\phi$. It also has a symmetrizing form $t_\phi$
defined as the specialization of the canonical form $t$.\\
\\
\begin{remark} \emph{Sometimes we describe the morphism $\phi$ by the
formula}
$$u_{\mathcal{C},j} \mapsto \zeta_{e_\mathcal{C}}^j x^{n_{\mathcal{C},j}}.$$
\emph{If now we set $q:=\zeta x$ for some root of unity $\zeta \in
\mu(K)$, then the cyclotomic specialization $\phi$ becomes a
$\zeta$-\emph{cyclotomic specialization} and $\mathcal{H}_\phi$ can
be also considered over $\mathbb{Z}_K[q,q^{-1}]$.}
\end{remark}
\begin{px}\label{spetsial}
\small{\emph{The \emph{spetsial} Hecke algebra $\mathcal{H}_q^s(W)$ is the
1-cyclotomic algebra obtained by the specialization
$$u_{\mathcal{C},0} \mapsto q,\,\, u_{\mathcal{C},j} \mapsto \zeta_{e_\mathcal{C}}^j \textrm{
for } 1 \leq j \leq e_\mathcal{C}-1, \textrm{ for all } \mathcal{C}
\in \mathcal{A}/W.$$ For example, if $W:=G_4$, then
$$\mathcal{H}_q^s(W)=<S,T \textrm{ }|\, STS=TST,
(S-q)(S^2+S+1)=(T-q)(T^2+T+1)=0>.$$}}
\end{px}

The following result is proved in \cite{Chlou} (remarks following Thm. 3.3.3):

\begin{proposition}\label{cyclotomic split semisimple}
The algebra $K(y)\mathcal{H}_\phi$ is split semisimple. 
\end{proposition}

When $y$ specializes to $1$, the algebra $K(y)\mathcal{H}_\phi$ specializes to the group
algebra $KW$ (the form $t_\phi$ becoming the canonical form on the
group algebra). Thus, by ``Tits' deformation theorem'', the
specialization $v_{\mathcal{C},j} \mapsto 1$ defines the following bijections
$$\begin{array}{ccccc}
    \textrm{Irr}(K(\textbf{v})\mathcal{H}) & \leftrightarrow & \textrm{Irr}(K(y)\mathcal{H}_\phi) & \leftrightarrow & \textrm{Irr}(W) \\
    \chi_{\textbf{v}} & \mapsto & \chi_{\phi} & \mapsto & \chi.
  \end{array}$$

The following result is an immediate consequence of Theorem
$\ref{Schur element generic}$.

\begin{proposition}\label{Schur element cyclotomic}
The Schur element $s_{\chi_\phi}(y)$ associated with the irreducible
character $\chi_\phi$ of $K(y)\mathcal{H}_\phi$ is a Laurent
polynomial in $y$ of the form
$$s_{\chi_\phi}(y)=\psi_{\chi,\phi} y^{a_{\chi,\phi}} \prod_{\Phi \in
C_K}\Phi(y)^{n_{\chi,\phi}}$$ where $\psi_{\chi,\phi} \in
\mathbb{Z}_K$, $a_{\chi,\phi} \in \mathbb{Z}$, $n_{\chi,\phi} \in
\mathbb{N}$ and $C_K$ is a set of $K$-cyclotomic polynomials.
\end{proposition}

\subsection{Rouquier blocks of the cyclotomic Hecke algebras}

\begin{definition}\label{Rouquier ring}
We call Rouquier ring of $K$ and denote by $\mathcal{R}_K(y)$ the
$\mathbb{Z}_K$-subalgebra of $K(y)$
$$\mathcal{R}_K(y):=\mathbb{Z}_K[y,y^{-1},(y^n-1)^{-1}_{n\geq 1}]$$
\end{definition}

Let $\phi: v_{\mathcal{C},j} \mapsto y^{n_{\mathcal{C},j}}$ be a
cyclotomic specialization and $\mathcal{H}_\phi$ the corresponding
cyclotomic Hecke algebra. The \emph{Rouquier blocks} of
$\mathcal{H}_\phi$ are the blocks of the algebra
$\mathcal{R}_K(y)\mathcal{H}_\phi$.\\
\\
\begin{remark}\emph{ It has been shown by Rouquier
\cite{Rou}, that if $W$ is a Weyl group and $\mathcal{H}_\phi$ is
obtained via the ``spetsial'' cyclotomic specialization (see example
$\ref{spetsial}$), then its Rouquier blocks coincide with the
``families of characters'' defined by Lusztig. Thus, the Rouquier
blocks play an essential role in the program ``Spets'' (see
\cite{BMM2}) whose ambition is to give to complex reflection groups
the role of Weyl groups of as yet mysterious structures.}
\end{remark}

The Rouquier ring has the following interesting properties. Their proof is given in
\cite{Chlou}, Proposition $3.4.2$.

\begin{proposition}\label{Some properties of the Rouquier
ring}\emph{(Some properties of the Rouquier ring)}
\begin{enumerate}
  \item The group of units $\mathcal{R}_K(y)^\times$ of the Rouquier ring $\mathcal{R}_K(y)$
  consists of the elements of the form
  $$u y^n \prod_{\Phi \in \mathrm{Cycl}(K)} \Phi(y)^{n_\phi},$$
  where $u \in \mathbb{Z}_K^\times$, $n, n_\phi \in \mathbb{Z}$, $\mathrm{Cycl}(K)$ is the set
  of $K$-cyclotomic polynomials and
  $n_\phi=0$ for all but a finite number of $\Phi$.
  \item The prime ideals of $\mathcal{R}_K(y)$ are
  \begin{itemize}
    \item the zero ideal $\{0\}$,
    \item the ideals of the form $\mathfrak{p}\mathcal{R}_K(y)$,
    where $\mathfrak{p}$ is a prime ideal of $\mathbb{Z}_K$,
    \item the ideals of the form $P(y)\mathcal{R}_K(y)$, where
    $P(y)$ is an irreducible element of $\mathbb{Z}_K[y]$ of degree
    at least $1$, prime to $y$ and to $\Phi(y)$ for all $\Phi \in
    \mathrm{Cycl}(K)$.
  \end{itemize}
  \item The Rouquier ring $\mathcal{R}_K(y)$ is a Dedekind ring.
\end{enumerate}
\end{proposition}

Now let us recall the form of the Schur elements of the cyclotomic
Hecke algebra $\mathcal{H}_\phi$ given in proposition $\ref{Schur
element cyclotomic}$. If $\chi_\phi$ is an irreducible character of
$K(y)\mathcal{H}_\phi$, then its Schur element $s_{\chi_\phi}(y)$ is
of the form
$$s_{\chi_\phi}(y)=\psi_{\chi,\phi} y^{a_{\chi,\phi}} \prod_{\Phi \in
C_K}\Phi(y)^{n_{\chi,\phi}}$$ where $\psi_{\chi,\phi} \in
\mathbb{Z}_K$, $a_{\chi,\phi} \in \mathbb{Z}$, $n_{\chi,\phi} \in
\mathbb{N}$ and $C_K$ is a set of $K$-cyclotomic polynomials.

\begin{definition}\label{bad}
A prime ideal $\mathfrak{p}$ of $\mathbb{Z}_K$ lying over a prime
number $p$ is $\phi$-bad for $W$, if there exists $\chi_\phi \in
\textrm{\emph{Irr}}(K(y)\mathcal{H}_\phi)$ with $\psi_{\chi,\phi}
\in \mathfrak{p}$. If $\mathfrak{p}$ is $\phi$-bad for $W$, we say
that $p$ is a $\phi$-bad prime number for $W$.
\end{definition}
\begin{remark}\emph{ If $W$ is a Weyl group and $\phi$ is the
``spetsial'' cyclotomic specialization, then the $\phi$-bad prime
ideals are the ideals generated by the bad prime numbers (in the
``usual'' sense) for $W$ (see \cite{GeRo}, 5.2).}
\end{remark}\

Note that if $\mathfrak{p}$ is $\phi$-bad for $W$, then $p$ must
divide the order of the group (since
$s_{\chi_\phi}(1)=|W|/\chi(1)$).\\

Let us denote by $\mathcal{O}$ the Rouquier ring. By proposition
$\ref{p-blocks}$, the Rouquier blocks of $\mathcal{H}_\phi$ are
unions of the blocks of $\mathcal{O}_\mathcal{P}\mathcal{H}_\phi$
for all prime ideals $\mathcal{P}$ of $\mathcal{O}$. However, in all
of the following cases, due to the form of the Schur elements, the
blocks of $\mathcal{O}_\mathcal{P}\mathcal{H}_\phi$ are singletons
(\ie $e_{\chi_\phi}=\chi_\phi^\vee /s_{\chi_\phi} \in
\mathcal{O}_\mathcal{P}\mathcal{H}_\phi$ for all $\chi_\phi \in
\mathrm{Irr}(K(y)\mathcal{H}_\phi)$):
\begin{itemize}
  \item $\mathcal{P}$ is the zero ideal $\{0\}$.
  \item $\mathcal{P}$ is of the form $P(y)\mathcal{O}$, where
$P(y)$ is an irreducible element of $\mathbb{Z}_K[y]$ of degree at
least $1$, prime to $y$ and to $\Phi(y)$ for all $\Phi \in
\mathrm{Cycl}(K)$.
  \item $\mathcal{P}$ is of the form $\mathfrak{p}\mathcal{O}$, where
$\mathfrak{p}$ is a prime ideal of $\mathbb{Z}_K$ which is not
$\phi$-bad for $W$.
\end{itemize}
Therefore, applying proposition $\ref{p-blocks}$, we obtain 

\begin{proposition}\label{Rouquier blocks and central characters}
Two characters $\chi,\psi \in \emph{Irr}(W)$ are in the same Rouquier block of $\mathcal{H}_\phi$
if and only if there exists a finite sequence
$\chi_0,\chi_1,\ldots,\chi_n \in \emph{Irr}(W)$ and a finite
sequence $\mathfrak{p}_1,\ldots,\mathfrak{p}_n$ of $\phi$-bad prime
ideals for $W$ such that
\begin{itemize}
  \item $\chi_0=\chi$ and $\chi_n=\psi$,
  \item for all $j$ $(1\leq j \leq n)$,\,\,
         the characters $\chi_{j-1}$ and $\chi_{j}$
         belong to the same block of
         $\mathcal{O}_{\mathfrak{p}_j\mathcal{O}}\mathcal{H}_\phi.$
 \end{itemize}                   
\end{proposition}

The above proposition implies that if we know the blocks of the algebra $\mathcal{O}_{\mathfrak{p}\mathcal{O}}\mathcal{H}_\phi$ for every $\phi$-bad prime ideal
$\mathfrak{p}$ for $W$,
then we know the Rouquier blocks of  $\mathcal{H}_\phi$. In order to determine the former, we
can use the following theorem (\cite{Chlou}, Thm.3.2.17):

\begin{theorem}\label{main theorem}
Let $A:=\mathbb{Z}_K[\textbf{\emph{v}},\textbf{\emph{v}}^{-1}]$ and $\mathfrak{p}$ be a prime
ideal of $\mathbb{Z}_K$. 
Let $M_1,\ldots,M_k$ be all the
$\mathfrak{p}$-essential monomials for $W$ such that $\phi(M_j)=1$
for all $j=1,\ldots,k$. Set $\mathfrak{q}_0:=\mathfrak{p}A$,
$\mathfrak{q}_j:=\mathfrak{p}A+(M_j-1)A$ for $j=1,\ldots,k$ and
$\mathcal{Q}:=\{\mathfrak{q}_0,\mathfrak{q}_1,\ldots,\mathfrak{q}_k\}$.
Two irreducible characters $\chi,\psi \in \textrm{\emph{Irr}}(W)$
are in the same block of $\mathcal{O}_{\mathfrak{p}\mathcal{O}}\mathcal{H}_\phi$ if
and only if there exist a finite sequence
$\chi_0,\chi_1,\ldots,\chi_n \in \textrm{\emph{Irr}}(W)$ and a
finite sequence $\mathfrak{q}_{j_1},\ldots,\mathfrak{q}_{j_n} \in
\mathcal{Q}$ such that
\begin{itemize}
  \item $\chi_0=\chi$ and $\chi_n=\psi$,
  \item for all $i$ $(1\leq i \leq n)$,  the characters $\chi_{i-1}$ and $\chi_i$ are
  in the same block of $A_{\mathfrak{q}_{j_i}}\mathcal{H}$.
\end{itemize}
\end{theorem}

Let $\mathfrak{p}$ be a prime ideal of $\mathbb{Z}_K$ and  $\phi: v_{\mathcal{C},j} \mapsto y^{n_{\mathcal{C},j}}$ a cyclotomic specialization.
If $M=\prod_{\mathcal{C},j}v_{\mathcal{C},j}^{a_{\mathcal{C},j}}$
is a $\mathfrak{p}$-essential monomial for $W$, then
$$\phi(M)=1 \Leftrightarrow \sum_{\mathcal{C},j}a_{\mathcal{C},j}n_{\mathcal{C},j}=0.$$
Set $m:=\sum_{\mathcal{C}\in \mathcal{A}/W}e_\mathcal{C}$. The
hyperplane defined in $\mathbb{C}^m$ by the relation
$$\sum_{\mathcal{C},j}a_{\mathcal{C},j}t_{\mathcal{C},j}=0,$$ where
$(t_ {\mathcal{C},j})_{ \mathcal{C},j}$ is a set of $m$
indeterminates, is called \emph{$\mathfrak{p}$-essential hyperplane}
for $W$. A hyperplane in $\mathbb{C}^m$ is called \emph{essential}
for $W$, if it is $\mathfrak{p}$-essential for some prime ideal
$\mathfrak{p}$ of $\mathbb{Z}_K$ (Respectively, a monomial is called \emph{essential}
for $W$, if it is $\mathfrak{p}$-essential for some prime ideal
$\mathfrak{p}$ of $\mathbb{Z}_K$).\\

Let $H$ be an essential hyperplane corresponding to the monomial $M$ and let $\mathfrak{p}$ be a prime ideal of $\mathbb{Z}_K$ .
We denote by $\mathcal{B}_\mathfrak{p}^H$ the partition of $\mathrm{Irr}(W)$ 
into blocks of $A_{\mathfrak{q}_M}\mathcal{H}$, where $\mathfrak{q}_M:=(M-1)A+\mathfrak{p}A$.
Moreover, we denote by $\mathcal{B}_\mathfrak{p}^\emptyset$ the partition of $\mathrm{Irr}(W)$ 
into blocks of $A_{\mathfrak{p}A}\mathcal{H}$.

\begin{definition}\label{rb associated with hyperplane}
Let $H$ be an essential hyperplane for $W$. We call Rouquier blocks associated with the hyperplane
$H$ (resp. with no essential hyperplane), and denote by  $\mathcal{B}^H$ (resp. $\mathcal{B}^\emptyset)$, the partition of $\mathrm{Irr}(W)$  generated  by the partitions
$\mathcal{B}_\mathfrak{p}^H$ (resp. $\mathcal{B}_\mathfrak{p}^\emptyset$), where $\mathfrak{p}$ runs over the set of prime ideals
of $\mathbb{Z}_K$.  
\end{definition}

With the help of the above definition and thanks to proposition $\ref{Rouquier blocks and central characters}$ and theorem $\ref{main theorem}$, we obtain the following characterization for the Rouquier blocks of a cyclotomic Hecke algebra:

\begin{proposition}\label{explain AllBlocks}
Let  $\phi: v_{\mathcal{C},j} \mapsto y^{n_{\mathcal{C},j}}$ be a cyclotomic specialization. The Rouquier blocks of the cyclotomic Hecke algebra $\mathcal{H}_\phi$ correspond to the partition of $\mathrm{Irr}(W)$  generated  by the partitions
$\mathcal{B}^H$, where $H$ runs over the set of all essential hyperplanes the integers
$n_{\mathcal{C},j}$ belong to. If the $n_{\mathcal{C},j}$ belong to no essential hyperplane, then 
the Rouquier blocks of $\mathcal{H}_\phi$ coincide with the partition  $\mathcal{B}^\emptyset$.
\end{proposition}

\begin{cordef}\label{assoc with hyp}Let
$\phi: v_{\mathcal{C},j} \mapsto y^{n_{\mathcal{C},j}}$ a cyclotomic specialization  such that  the integers
$n_{\mathcal{C},j}$ belong to only one essential hyperplane $H$ (resp. to no essential hyperplane). We say that $\phi$ is a cyclotomic specialization associated with the essential hyperplane $H$ (resp. with no essential hyperplane). Then the Rouquier blocks of $\mathcal{H}_\phi$ coincide with the partition   $\mathcal{B}^H$ (resp.  $\mathcal{B}^\emptyset$). 
\end{cordef}

By taking cyclotomic specializations associated to each (or no) essential hyperplane and calculating the Rouquier blocks of the corresponding cyclotomic Hecke algebras, we have been able to determine the Rouquier blocks for all exceptional complex reflection groups in \cite{Chlou}. In the next section, we are going to do the same for the group $G(d,1,r)$.

 \subsection{Functions $a$ and $A$}

Following the notations in \cite{BMM2}, 6B, for every element $P(y)
\in \mathbb{C}(y)$, we call
\begin{itemize}
  \item \emph{valuation of $P(y)$ at $y$} and denote by $\mathrm{val}_y(P)$ the order of $P(y)$
  at 0 (we have $\mathrm{val}_y(P)<0$ if 0 is a pole of $P(y)$ and $\mathrm{val}_y(P)>0$ if 0 is a zero of $P(y)$),
  \item \emph{degree of $P(y)$ at $y$} and denote by $\mathrm{deg}_y(P)$ the opposite of the
  valuation of $P(1/y)$.
\end{itemize}
Moreover, if $x:=y^{|\mu(K)|}$, then
$$\mathrm{val}_x(P):=\frac{\mathrm{val}_y(P)}{|\mu(K)|} \textrm{ and }
\mathrm{deg}_x(P):=\frac{\mathrm{deg}_y(P)}{|\mu(K)|}.$$ 
For $\chi
\in \mathrm{Irr}(W)$, we define
$$a_{\chi_\phi}:=\mathrm{val}_x(s_{\chi_\phi}(y)) \,\textrm{ and }\,
A_{\chi_\phi}:=\mathrm{deg}_x(s_{\chi_\phi}(y)).$$ 
The following
result is proven in \cite{BK}, Prop.2.9.

\begin{proposition}\label{aA}\
Let $\chi,\psi \in \mathrm{Irr}(W)$. If $\chi_\phi$ and
        $\psi_\phi$ belong to the same Rouquier block, then
        $$a_{\chi_\phi}+A_{\chi_\phi}=a_{\psi_\phi}+A_{\psi_\phi}.$$
\end{proposition}

\section{Rouquier blocks for the Ariki-Koike algebras}

We will start this section by introducing some notations and results in combinatorics (cf. \cite{BK}, \S 3A) which will be useful for the description of the Rouquier blocks of the cyclotomic Ariki-Koike algebras, \ie the cyclotomic Hecke algebras associated to the group $G(d,1,r)$.

\subsection{Combinatorics}

Let $\el=(\el_1,\el_2,\ldots,\el_h)$ be a partition, \ie a finite decreasing sequence of positive integers:
$$\el_1 \geq \el_2 \geq \ldots \geq \el_h \geq 1.$$
The integer
$$|\el|:=\el_1+\el_2+\ldots+\el_h$$
is called \emph{the size of $\el$}. We also say that $\lambda$ \emph{is a partition of }
$|\el|$.
The integer $h$ is called \emph{the height of $\el$} and we set $h_\el:=h$. To each partition $\el$ we associate its \emph{$\beta$-number}, $\eb_\el=(\eb_1,\eb_2,\ldots,\eb_h)$, defined by
$$\eb_1:=h+\el_1-1,\eb_2:=h+\el_2-2,\ldots,\eb_h:=h+\el_h-h.$$

\subsection*{\normalsize Multipartitions}

From now on,  $d$ is a positive integer. Let $\el=(\el^{(0)},\el^{(1)},\ldots,\el^{(d-1)})$ be a $d$-partition, \ie a family of $d$ partitions indexed by the set $\{0,1,\ldots,d-1\}$. We set 
$$h^{(a)}:=h_{\el^{(a)}}, \,\,\, \eb^{(a)}:=\eb_{\el^{(a)}}$$
and we have
$$ \el^{(a)}=(\el_1^{(a)},\el_2^{(a)},\ldots,\el_{h^{(a)}}^{(a)}).$$
The integer
$$|\el|:=\sum_{a=0}^{d-1}|\el^{(a)}|$$
is called \emph{the size of $\el$}. We also say that $\lambda$ \emph{is a $d$-partition of}
$|\el|$.

\subsection*{\normalsize Ordinary symbols}

If $\eb=(\eb_1,\eb_2,\ldots,\eb_h)$ is a sequence of positive integers such that $\eb_1>\eb_2>\ldots>\eb_h$ and $m$ is a positive integer, then the $m$-``shifted'' of $\eb$ is the sequence of numbers defined by
$$\eb[m]=(\eb_1+m,\eb_2+m,\ldots,\eb_h+m,m-1,m-2,\ldots,1,0).$$

Let  $\el=(\el^{(0)},\el^{(1)},\ldots,\el^{(d-1)})$ be a $d$-partition.  We call \emph{$d$-height of $\el$} the family $(h^{(0)},h^{(1)},\ldots,h^{(d-1)})$ and we define the 
 \emph{height of $\el$} to be the integer
$$h_\el:=\mathrm{max}\,\{h^{(a)} \,|\, (0 \leq a \leq d-1)\}.$$

\begin{definition}\label{ordinary standard symbol}
The ordinary standard symbol of $\el$ is the family of numbers defined by
$$B_\el=(B_\el^{(0)},B_\el^{(1)},\ldots,B_\el^{(d-1)}),$$
where, for all $a$ $(0 \leq a \leq d-1)$, we have
$$B_\el^{(a)}:=\eb^{(a)}[h_\el-h^{(a)}].$$
An ordinary symbol of $\el$ is a symbol obtained from the ordinary standard symbol by shifting all the rows by the same integer.
\end{definition}

The ordinary standard symbol of a $d$-partition $\el$ is of the form
$$
\begin{array}{cccccc}
 B_\el^{(0)} & = & b_1^{(0)} & b_2^{(0)} & \ldots & b_{h_\el}^{(0)} \\
 B_\el^{(1)} & = & b_1^{(1)} & b_2^{(1)} & \ldots & b_{h_\el}^{(1)} \\
  \vdots  &    & \vdots       &\vdots         &\vdots & \vdots  \\
 B_\el^{(d-1)} & = & b_1^{(d-1)} & b_2^{(d-1)} & \ldots & b_{h_\el}^{(d-1)} \\
\end{array}
$$

The \emph{ordinary content} of a $d$-partition of ordinary standard symbol $B$ is the ``set with repetition''
$$\mathrm{Cont}_\el = B_\el^{(0)} \cup B_\el^{(1)} \cup \ldots \cup B_\el^{(d-1)}$$
or (with the above notations) the polynomial defined by
$$\mathrm{Cont}_\el(x):= \sum_{a,i} x^{b_i^{(a)}}.$$

\begin{px}{\small \emph{Let us take $d=2$ and $\el=((2,1),(3))$. Then}
$$
B_\el=
\left(
\begin{array}{cc}
  3 & 1   \\
  4 & 0    
\end{array}
\right)
$$
\emph{\small We have
$\mathrm{Cont}_\el=\{0,1,3,4\}$ or $\mathrm{Cont}_\el(x)=1+x+x^3+x^4.$}}
\end{px}

\subsection*{\normalsize Charged symbols}

From now on, we suppose that we have a given ``weight system'', \ie a family of integers
$$m:=(m^{(0)},m^{(1)},\ldots,m^{(d-1)}).$$

Let  $\el=(\el^{(0)},\el^{(1)},\ldots,\el^{(d-1)})$ be a $d$-partition.  We call \emph{$(d,m)$-charged height of $\el$} the family $(hc^{(0)},hc^{(1)},\ldots,hc^{(d-1)})$, where
$$hc^{(0)}:=h^{(0)}-m^{(0)},hc^{(1)}:=h^{(1)}-m^{(1)},\ldots,hc^{(d-1)}:=h^{(d-1)}-m^{(d-1)}.$$
We define the 
 \emph{$m$-charged height of $\el$} to be the integer
$$hc_\el:=\mathrm{max}\,\{hc^{(a)} \,|\, (0 \leq a \leq d-1)\}.$$

\begin{definition}\label{charged standard symbol}
The $m$-charged standard symbol of $\el$ is the family of numbers defined by
$$Bc_\el=(Bc_\el^{(0)},Bc_\el^{(1)},\ldots,Bc_\el^{(d-1)}),$$
where, for all $a$ $(0 \leq a \leq d-1)$, we have
$$Bc_\el^{(a)}:=\eb^{(a)}[hc_\el-hc^{(a)}].$$
An $m$-charged symbol of $\el$ is a symbol obtained from the $m$-charged standard symbol by shifting all the rows by the same integer.
\end{definition}
\begin{remark} \emph{The ordinary symbols correspond to the weight system }
\begin{center}
$m^{(0)}=m^{(1)}=\ldots=m^{(d-1)}=0.$
\end{center}
\end{remark}

The $m$-charged standard symbol of $\el$ is a tableau of numbers arranged into $d$ rows indexed by the set $\{0,1,\ldots,d-1\}$ such that the $a^{\mathrm{th}}$ row has length equal to $hc_\el+m^{(a)}$. For all $a$ $(0 \leq a \leq d-1)$, we set $l^{(a)}:=hc_\el+m^{(a)}$ and we denote by
$$
\begin{array}{cccccc}
 Bc_\el^{(a)} & = & bc_1^{(a)} & bc_2^{(a)} & \ldots & bc_{l^{(a)}}^{(a)} 
 \end{array}$$
the $a^{\mathrm{th}}$ row of the $m$-charged standard symbol.

The \emph{$m$-charged content} of a $d$-partition of $m$-charged standard symbol $Bc$ is the ``set with repetition''
$$\mathrm{Contc}_\el = Bc_\el^{(0)} \cup Bc_\el^{(1)} \cup \ldots \cup Bc_\el^{(d-1)}$$
or (with the above notations) the polynomial defined by
$$\mathrm{Contc}_\el(x):= \sum_{a,i} x^{bc_i^{(a)}}.$$

\begin{px}{\small \emph{Let us take $d=2$, $\el=((2,1),(3))$ and $m=(-1,2)$. Then}
$$
Bc_\el=
\left(
\begin{array}{ccccc}
  3 & 1 &     &    & \\
  7 & 3 & 2 & 1 & 0    
\end{array}
\right)
$$
\emph{\small We have
$\mathrm{Contc}_\el=\{0,1,1,2,3,3,7\}$ or $\mathrm{Contc}_\el(x)=1+2x+x^2+2x^3+x^7.$}}
\end{px}

\subsection{Generic Ariki-Koike algebras}

The group $G(d,1,r)$ is the group of all monomial $r \times r$ matrices with entries in $\mu_d$. It is isomorphic to the wreath product $\mu_d \wr
\mathfrak{S}_r$  and its field of definition is $\mathbb{Q}(\zeta_d)$.

The \emph{generic Ariki-Koike algebra} of $G(d,1,r)$ (cf. \cite{ArKo}, \cite{BM}) is the algebra $\mathcal{H}_{d,r}$ generated over the Laurent ring of polynomials in $d+1$ indeterminates  
$$\mathcal{O}_d:=\mathbb{Z}[u_0,u_0^{-1},u_1,u_1^{-1},\ldots,u_{d-1},u_{d-1}^{-1},x,x^{-1}]$$
by the elements $\mathrm{\textbf{s}},\mathrm{\textbf{t}}_1,\mathrm{\textbf{t}}_2,\ldots,\mathrm{\textbf{t}}_{r-1}$ satisfying the relations
\begin{itemize}
\item $\mathrm{\textbf{s}}\mathrm{\textbf{t}}_1\mathrm{\textbf{s}}\mathrm{\textbf{t}}_1=\mathrm{\textbf{t}}_1\mathrm{\textbf{s}}\mathrm{\textbf{t}}_1\mathrm{\textbf{s}}$, $\mathrm{\textbf{s}}\mathrm{\textbf{t}}_j=\mathrm{\textbf{t}}_j\mathrm{\textbf{s}} \textrm{ for } j\neq 1$,
\item $\mathrm{\textbf{t}}_j\mathrm{\textbf{t}}_{j+1}\mathrm{\textbf{t}}_j=\mathrm{\textbf{t}}_{j+1}\mathrm{\textbf{t}}_j\mathrm{\textbf{t}}_{j+1}$,  $ \mathrm{\textbf{t}}_i\mathrm{\textbf{t}}_j=\mathrm{\textbf{t}}_j\mathrm{\textbf{t}}_i \textrm{ for } |i-j|>1$,
\item $(\mathrm{\textbf{s}}-u_0)(\mathrm{\textbf{s}}-u_1)\ldots(\mathrm{\textbf{s}}-u_{d-1})=(\mathrm{\textbf{t}}_j-x)(\mathrm{\textbf{t}}_j+1)=0$.
\end{itemize}\

For every $d$-partition $\el=(\el^{(0)}, \el^{(1)}, \ldots, \el^{(d-1)})$ of $r$, we consider the free $\mathcal{O}_d$-module which has as basis the family of standard tableaux of $\el$. We can give to this module the structure of a $\mathcal{H}_{d,r}$-module (cf. \cite{ArKo}, \cite{Ar1}, \cite{GrLe}) and thus obtain the \emph{Specht module} \textbf{Sp}$^\el$ associated to $\el$.  

Set $\mathcal{K}_d:=\mathbb{Q}(u_0,u_1,\ldots,u_{d-1},x)$ the field of fractions of $\mathcal{O}_d$.
The  $\mathcal{K}_d\mathcal{H}_{d,r}$-module  $\mathcal{K}_d$\textbf{Sp}$^\el$, obtained by extension of scalars, is absolutely irreducible and every irreducible $\mathcal{K}_d\mathcal{H}_{d,r}$-module is isomorphic to a module of this type. Thus $\mathcal{K}_d$ is a splitting field for $\mathcal{H}_{d,r}$. 
We denote by $\chi_\el$ the (absolutely) irreducible character of the $\mathcal{K}_d\mathcal{H}_{d,r}$-module  \textbf{Sp}$^\el$.\\

Since the algebra $\mathcal{K}_d\mathcal{H}_{d,r}$ is split semisimple, the Schur elements of its irreducible characters belong to $\mathcal{O}_d$. The following result by Mathas (\cite{Mat}, Cor.6.5) gives a description of the Schur elements. The same result has been obtained independently by Geck, Iancu and Malle in \cite{GIM}.

\begin{proposition}\label{schur}
Let $\el$ be a $d$-partition of $r$ with ordinary standard symbol
$B_\el=(B_\el^{(0)}, B_\el^{(1)}, \ldots, B_\el^{(d-1)})$. Fix $L \geq h_\el$, where $h_\el$ is the height of $\el$.
We set $B_{\el,L}:=(B_\el^{(0)}[L-h_\el], B_\el^{(1)}[L-h_\el], \ldots, B_\el^{(d-1)}[L-h_\el])=(B_{\el,L}^{(0)}, B^{(1)}_{\el,L}, \ldots, B^{(d-1)}_{\el,L})$ and
$B_{\el,L}^{(s)}=(b_1^{(s)},b_2^{(s)},\ldots,
b_{L}^{(s)})$.  Let $a_{L}:=r(d-1)+\binom{ d}{ 2}\binom{ L}{ 2}$ and $b_{L}:=dL(L-1)(2dL-d-3)/12$. Then the Schur element of the irreducible character $\chi_\el$ is given by the formulae $s_\el=(-1)^{a_{L}} x^{b_{L}}(x-1)^{-r}(u_0u_1\ldots u_{d-1})^{-r}\nu_\el/ \delta_\el$, where
$$
\nu_\el=
\prod_{0\leq s<t<d}(u_s-u_t)^L\prod_{0 \leq s,t <d}\prod_{b_s \in B_{\el,L}^{(s)}}\prod_{1 \leq k \leq b_s} 
(x^ku_s-u_t)$$
and
$$\delta_\el=\prod_{0\leq s< t <d}\prod_{(b_s,b_t) \in B_{\el,L}^{(s)}\times B_{\el,L}^{(t)}}(x^{b_s}u_s-x^{b_t}u_t) \prod_{0 \leq s <d} \prod_{1 \leq i < j \leq L}(x^{b_i^{(s)}}u_s-x^{b_j^{(s)}}u_s).
$$
\end{proposition}

We have already mentioned that the field of definition of $G(d,1,r)$ is $K:=\mathbb{Q}(\zeta_d)$. If we set
$$v_j^{|\mu(K)|}:=\zeta_d^{-j}u_j (0 \leq j \leq d-1) \textrm{ and } z^{|\mu(K)|}:=x,$$
then theorem $\ref{Semisimplicity Malle}$ implies that the algebra
$K(v_0,v_1,\ldots,v_{d-1},z)\mathcal{H}_{d,r}$ is split semisimple.
Proposition $\ref{schur}$ implies that the essential monomials for  $G(d,1,r)$ are of the form 
\begin{itemize}
\item $z^kv_sv_t^{-1}$ for $0 \leq s<t<d$ and $-r<k <r$, 
\item $z$
\end{itemize}
\begin{remark}\emph{ The monomial $z$ can be seen as a monomial of the form $z_0z_1^{-1}$, if, in the definition of the Ariki-Koike algebra, we replace the relation}
$$({\textbf{t}}_j-x)(\mathrm{\textbf{t}}_j+1)=0
\textrm{ \emph{by} }
({\textbf{t}}_j-x_0)(\mathrm{\textbf{t}}_j+x_1)=0$$
\emph{and we set}
$$z_0^{|\mu(K)|}:=x_0 \textrm{ \emph{and} } z_1^{|\mu(K)|}:=x_1.$$
\end{remark} 

\subsection{Cyclotomic Ariki-Koike algebras}

Let $y$ be an indeterminate and let
$$\phi : \left\{ 
\begin{array}{ll} 
v_j \mapsto y^{m_j}, (0 \leq j <d),\\ 
z \mapsto y^n
\end{array} \right. 
$$
be a cyclotomic specialization. If we set $q:=y^{|\mu(K)|}$,
then $\phi$ can be described as follows
$$\phi : \left\{ 
\begin{array}{ll} 
u_j \mapsto \zeta_d^j q^{m_j}, (0 \leq j <d),\\ 
x \mapsto q^n
\end{array} \right. 
$$

The corresponding cyclotomic Hecke algebra $(\mathcal{H}_{d,r})_\phi$ can be considered either over the ring $\mathbb{Z}_{K}[y,y^{-1}]$ or over the ring $\mathbb{Z}_{K}[q,q^{-1}]$. We define the Rouquier blocks of  $(\mathcal{H}_{d,r})_\phi$ to be the blocks of  $(\mathcal{H}_{d,r})_\phi$ defined over the Rouquier ring $\mathcal{R}_{K}(y)$ in $K(y)$. However, in other texts, as, for example, in \cite{BK}, the Rouquier blocks are determined over the Rouquier ring $\mathcal{R}_{K}(q)$ in $K(q)$. Since $\mathcal{R}_{K}(y)$ is the integral closure of $\mathcal{R}_{K}(q)$ in the splitting field  $K(y)$ for $(\mathcal{H}_{d,r})_\phi$, proposition $\ref{Galois action on integral closure}$ establishes a relation between the blocks of
$\mathcal{R}_{K}(y)(\mathcal{H}_{d,r})_\phi$ and the blocks of $\mathcal{R}_{K}(q)(\mathcal{H}_{d,r})_\phi$. Moreover, in our case we can prove that

\begin{proposition}\label{q or y}
The  blocks of
$\mathcal{R}_{K}(y)(\mathcal{H}_{d,r})_\phi$ and the blocks of $\mathcal{R}_{K}(q)(\mathcal{H}_{d,r})_\phi$ coincide.
\end{proposition}
\begin{apod}{By proposition $\ref{Galois action on integral closure}$, we know that the blocks of  $\mathcal{R}_{K}(q)(\mathcal{H}_{d,r})_\phi$ are unions of the blocks of $\mathcal{R}_{K}(y)(\mathcal{H}_{d,r})_\phi$. Now let $e$ be a block-idempotent of $\mathcal{R}_{K}(y)(\mathcal{H}_{d,r})_\phi$. Since $\mathcal{K}_d$ is a splitting field for $\mathcal{H}_{d,r}$, proposition $\ref{schur elements and idempotents}$ implies that $e$ belongs to $K(q)(\mathcal{H}_{d,r})_\phi$. Thus we have
$$e \in  \mathcal{R}_{K}(y)(\mathcal{H}_{d,r})_\phi \cap K(q)(\mathcal{H}_{d,r})_\phi
=\mathcal{R}_{K}(q)(\mathcal{H}_{d,r})_\phi,$$
since the ring $\mathcal{R}_{K}(q)$ is integrally closed and $\mathcal{R}_{K}(y)$ is integral over it
($y^{|\mu(K)|}-q=0$). Thus, $e$ is a sum of blocks of $\mathcal{R}_K(q)(\mathcal{H}_{d,r})_\phi$.}
\end{apod}

\subsection{Residue equivalency}

Let $\phi$ be a cyclotomic specialization like above and set  $\mathcal{O}:=\mathcal{R}_K(q)$. Following proposition $\ref{Rouquier blocks and central characters}$, in order to obtain the Rouquier blocks of $(\mathcal{H}_{d,r})_\phi$, we need to calculate the blocks of
$\mathcal{O}_{\mathfrak{p}\mathcal{O}}(\mathcal{H}_{d,r})_\phi$ for all
$\phi$-bad prime ideals $\mathfrak{p}$ of $\mathbb{Z}_K$.

Let $\mathfrak{p}$ be a prime ideal of  $\mathbb{Z}_K$ lying over a prime number $p$. By proposition $\ref{Some properties of the Rouquier ring}$ the ring $\mathcal{O}$ is a Dedekind ring and thus $\mathcal{O}_{\mathfrak{p}\mathcal{O}}$ is a discrete valuation ring. If we denote by $k_\mathfrak{p}$ its residue field, then, by proposition $\ref{residue field}$, the blocks of
$\mathcal{O}_{\mathfrak{p}\mathcal{O}}(\mathcal{H}_{d,r})_\phi$ coincide with the blocks of
$k_{\mathfrak{p}}(\mathcal{H}_{d,r})_\phi$. 
We denote by $\pi_\mathfrak{p}$ the surjective map $\mathcal{O}_{\mathfrak{p}\mathcal{O}}\twoheadrightarrow k_\mathfrak{p}$.

\begin{definition}\label{diagram of a multipartition}
The diagram of a $d$-partition $\lambda$ is the set
$$[\el]:=\{(i,j,a)\,|\,  (0 \leq a \leq d-1)(1\leq i \leq h^{(a)})(1 \leq j \leq \el_i^{(a)})\}.$$
A node is any ordered triple $(i,j,a)$. 
\end{definition}

The $\mathfrak{p}$-\emph{residue} of the node $x=(i,j,a)$ with respect to  $\phi$ is

$$\mathrm{res}_{\mathfrak{p},\phi}(x)= \left\{ 
\begin{array}{ll} 
\pi_\mathfrak{p}(\zeta_d^a q^{m_a}q^{n(j-i)} ) & \text{if }  n \neq 0,\\  \\
(\pi_\mathfrak{p}(j-i), \zeta_d^a q^{m_a}) & \text{if } n=0 \text{ and } \pi_\mathfrak{p}(\zeta_d^a q^{m_a}) \neq \pi_\mathfrak{p}(\zeta_d^b q^{m_b})  \text{ for } a\neq b,\\  \\
\pi_\mathfrak{p}(\zeta_d^a q^{m_a}) & \text{otherwise. }  
\end{array} \right. 
$$

Let $\mathrm{Res}_{\mathfrak{p},\phi}:=\{\mathrm{res}_{\mathfrak{p},\phi}(x)\,|\,x \in [\el] \textrm{ for some $d$-partition $\el$ of $r$}\}$ be the set of all possible residues. For any $d$-partition $\el$ of $r$ and $f \in \mathrm{Res}_{\mathfrak{p},\phi}$, we define
$$C_f(\el)=\# \{x \in [\el] \,|\, \mathrm{res}(x)=f\}.$$

Adapting definition $2.10$ of \cite{LyMa}, we obtain

\begin{definition}\label{residue equivalent}
Let $\el$ and $\mu$ be two $d$-partitions of $r$. We say that $\el$ and $\mu$ are $\mathfrak{p}$-residue equivalent  with respect to $\phi$ if $C_f(\el)=C_f(\mu)$ for all $f \in \mathrm{Res}_{\mathfrak{p},\phi}$.
\end{definition}

Then \cite{LyMa}, Theorem 2.13 implies that

\begin{theorem}\label{thm residue}
Two irreducible characters $(\chi_\el)_\phi$ and $(\chi_\mu)_\phi$ are in the same block of $k_\mathfrak{p}(\mathcal{H}_{d,r})_\phi$ if and only if $\el$ and $\mu$ are $\mathfrak{p}$-residue equivalent with respect to $\phi$.
\end{theorem}

The above result, in combination with proposition $\ref{Rouquier blocks and central characters}$ gives

\begin{corollary}\label{thm A}
Two irreducible characters $(\chi_\el)_\phi$ and $(\chi_\mu)_\phi$ are in the same Rouquier block of 
$(\mathcal{H}_{d,r})_\phi$
if and only if there exists a finite sequence
$\el_{(0)},\el_{(1)},\ldots,\el_{(m)}$ of $d$-partitions of $r$ and a finite
sequence $\mathfrak{p}_1,\ldots,\mathfrak{p}_m$ of $\phi$-bad prime
ideals for $W$ such that
\begin{itemize}
  \item $\el_{(0)}=\el$ and $\el_{(m)}=\mu$,
  \item for all $j$ $(1\leq j \leq m)$,\,\,
         the $d$-partitions $\el_{(j-1)}$ and $\el_{(j)}$
         are $\mathfrak{p}_j$-residue equivalent with respect to $\phi$.
 \end{itemize}                   
\end{corollary}

\subsection{Rouquier blocks and charged content}

Theorem $3.13$ in \cite{BK} gives a description of the Rouquier blocks of the cyclotomic Ariki-Koike algebras when $n \neq 0$. However, in the proof it is supposed that $1-\zeta_d$ always belongs to a prime ideal of $\mathbb{Z}[\zeta_d]$. This is not correct, unless $d$ is the power of a prime number. Therefore, we will state here the part of the theorem that is correct and only for the case $n=1$.

\begin{theorem}\label{BK}
Let $\phi$ be a cyclotomic specialization such that $\phi(x) = q$. 
If two irreducible characters $(\chi_\el)_\phi$ and $(\chi_\mu)_\phi$
are in the same Rouquier block of $(\mathcal{H}_{d,r})_\phi$, then $\mathrm{Contc}_\el=\mathrm{Contc}_\mu$ with respect to the weight system $m=(m_0,m_1,\ldots,m_{d-1})$. The converse holds when $d$ is the power of a prime number.
\end{theorem}

\subsection{Determination of the Rouquier blocks}

In this section, we are going to determine the Rouquier blocks for all cyclotomic Ariki-Koike algebras by determining the Rouquier blocks associated with no and each essential hyperplane for $G(d,1,r)$.
Due to corollary $\ref{assoc with hyp}$, it suffices to consider a cyclotomic specialization associated with no and each essential hyperplane and calculate the Rouquier blocks of the corresponding cyclotomic Hecke algebra. Following the description of the essential monomials in section $3.2$, we obtain that the essential hyperplanes for $G(d,1,r)$ are of the form
\begin{itemize}
\item $kN+M_s-M_t=0$ for $0 \leq s<t<d$ and $-r<k<r$.
\item $N=0$.
\end{itemize} 

\subsection*{\normalsize Case 1: No essential hyperplane}

If $\phi$ is a cyclotomic specialization associated with no essential hyperplane, then the description of the Schur elements by proposition $\ref{schur}$ implies that there are no $\phi$-bad prime ideals for $G(d,1,r)$. Therefore, every irreducible character is a block by itself.

\begin{proposition}\label{no essential hyperplane}
The Rouquier blocks associated with no essential hyperplane are trivial.
\end{proposition}

\subsection*{\normalsize Case 2: Essential hyperplane of the form $kN+M_s-M_t=0$}

The following result is an immediate consequence of the description of the Schur elements by proposition $\ref{schur}$.

\begin{proposition}\label{existence of prime ideal}
Let $s,t,k$ be three integers such that $0 \leq s<t<d$ and $-r<k<r$. The hyperplane
$$H:\,kN+M_s-M_t=0$$
is essential for $G(d,1,r)$ if and only if there exists a prime ideal $\mathfrak{p}$ of $\mathbb{Z}[\zeta_d]$ such that $\zeta_d^s-\zeta_d^t \in \mathfrak{p}$. Moreover, in this case, $H$ is $\mathfrak{p}$-essential for $G(d,1,r)$.
\end{proposition}

\begin{px}\emph{\small The hyperplane $M_0=M_1$ is $2$-essential for $G(2,1,r)$, whereas it isn't essential for $G(6,1,r)$, for all $r>0$.}
\end{px}

From now on, we assume that $kN+M_s-M_t=0$ is an essential hyperplane for $G(d,1,r)$, \ie that there exists a prime ideal $\mathfrak{p}$ of $\mathbb{Z}[\zeta_d]$ such that $\zeta_d^s-\zeta_d^t \in \mathfrak{p}$. Let $$\phi : \left\{ 
\begin{array}{ll} 
u_j \mapsto \zeta_d^j q^{m_j}, (0 \leq j <d),\\ 
x \mapsto q^n
\end{array} \right. 
$$
 be a cyclotomic specialization associated with that essential hyperplane. Our aim is the determination of the Rouquier blocks of $(\mathcal{H}_{d,r})_{\phi}$.

For the notations used in the following theorem, the reader should refer to section $3.1$.

\begin{proposition}\label{essential hyperplane of type 1}
Let $\el$, $\mu$ be two $d$-partitions of $r$. The irreducible characters $(\chi_\el)_{\phi}$ and $(\chi_\mu)_{\phi}$ are in the same Rouquier block of $(\mathcal{H}_{d,r})_{\phi}$ if and only if the following conditions are satisfied:
\begin{enumerate}
\item We have $\el^{(a)}=\mu^{(a)}$ for all $a \notin \{s,t\}.$
\item If $\el^{st}:=(\el^{(s)},\el^{(t)})$ and $\mu^{st}:=(\mu^{(s)},\mu^{(t)})$, then
$\mathrm{Contc}_{\el^{st}}= \mathrm{Contc}_{\mu^{st}}$ with respect to the weight system $(0,k)$.
\end{enumerate}
\end{proposition}
\begin{apod}{We can assume, without loss of generality, that $n=1$. We can also assume that $m_s=0$ and $m_t=k$.

Suppose that $(\chi_\el)_{\phi}$ and $(\chi_\mu)_{\phi}$ belong to  the same Rouquier block of $(\mathcal{H}_{d,r})_{\phi}$. Due to theorem $\ref{BK}$, we have $\mathrm{Contc}_\el=\mathrm{Contc}_\mu$ with respect to the weight system $m=(m_0,m_1,\ldots,m_{d-1})$. Since the $m_a$, $a \notin \{s,t\}$ could take any value (as long as they don't belong to another essential hyperplane), we must have that 
$\el^{(a)}=\mu^{(a)}$ for all $a \notin \{s,t\}.$ Moreover, the equality $\mathrm{Contc}_\el=\mathrm{Contc}_\mu$ implies that the corresponding $m$-charged standard symbols $Bc_\el$ and $Bc_\mu$ have the same cardinality and thus $hc_\el=hc_\mu$. Therefore, we obtain
$$Bc_\el^{(a)}=\beta^{(a)}_\el[hc_\el-hc_\el^{(a)}]=\beta^{(a)}_\mu [hc_\mu-hc_\mu^{(a)}]=Bc_\mu^{(a)} \textrm{ for all } a \notin \{s,t\}.$$
Consequently, we have the following equality between ``sets with repetition'':
$$Bc_\el^{(s)} \cup Bc_\el^{(t)}=Bc_\mu^{(s)} \cup Bc_\mu^{(t)}.$$
We can assume that the $m_a$, $a \notin \{s,t\}$ are sufficiently large so that 
$$hc_\el \in \{hc_\el^{(s)},hc_\el^{(t)}\} \textrm{ and } hc_\mu \in \{hc_\mu^{(s)},hc_\mu^{(t)}\}.$$
In this case, if  $\el^{st}:=(\el^{(s)},\el^{(t)})$ and $\mu^{st}:=(\mu^{(s)},\mu^{(t)})$, then
$$Bc_{\el^{st}}^{(0)}=Bc_\el^{(s)},  Bc_{\el^{st}}^{(1)}=Bc_\el^{(t)}, 
Bc_{\mu^{st}}^{(0)}=Bc_\mu^{(s)}, Bc_{\mu^{st}}^{(1)}=Bc_\mu^{(t)}$$
and we obtain $\mathrm{Contc}_{\el^{st}}= \mathrm{Contc}_{\mu^{st}}$ with respect to the weight system $(m_s,m_t)$. 

Now let us suppose that the conditions $1$ and $2$ are satisfied. Set $l:=|\el^{st}|$. Due to the first condition, we have $|\mu^{st}|=l$.   Let $\mathcal{H}_{2,l}$ be the generic Ariki-Koike algebra of the group $G(2,1,l)$ defined over the ring $$\mathbb{Z}[U_s,U_s^{-1},U_t,U_t^{-1},X,X^{-1}].$$ The group $G(2,1,l)$ is isomorphic to the cyclic group of order $2$ for $l=1$ and to the Coxeter group $B_l$ for $l \geq 2$.
Let us consider the cyclotomic specialization $$\vartheta:U_s\mapsto q^{m_s}, U_t \mapsto -q^{m_t}, X \mapsto q.$$

Due to theorem $\ref{BK}$, the condition $2$ implies that the characters
 $(\chi_{\lambda^{st}})_{\vartheta}$ and  $(\chi_{\mu^{st}})_{\vartheta}$ belong to the same Rouquier block of $(\mathcal{H}_{2,l})_{\vartheta}$. Therefore, we obtain that $kN+M_s-M_t=0$ is a $2$-essential hyperplane for $G(2,1,l)$ and that, due to corollary $\ref{thm A}$, 
 $\el^{st}$ and $\mu^{st}$ are $2$-residue equivalent with respect to $\vartheta$. 
In order to check whether $\el$ and $\mu$ are $\mathfrak{p}$-residue equivalent with respect to $\phi$, we only need to consider the nodes with third entry $s$ or $t$ (thanks to condition $1$). The nodes of $\el$ (resp. of $\mu$) with third entry $s$ or $t$ are the nodes of $\el^{st}$ (resp. $\mu^{st}$).
The $\mathfrak{p}$-residues of these nodes with respect to $\phi$ can be obtained by replacing  $q^{m_s}$ by $\zeta_d^s q^{m_s}$ and  $-q^{m_t}$ by $\zeta_d^t q^{m_t}$ into the $2$-residues with respect to $\vartheta$ of the nodes belonging to $[\el^{st}]$ and $[\mu^{st}]$. Since
$\el^{st}$ and $\mu^{st}$ are $2$-residue equivalent and $\zeta_d^s-\zeta_d^t \in \mathfrak{p}$ (when before we had $1-(-1) \in (2)$),  we obtain that  $\el$ and $\mu$ are $\mathfrak{p}$-residue equivalent with respect to $\phi$. Thus, by corollary $\ref{thm A}$, 
$(\chi_\el)_{\phi}$ and
$(\chi_{\mu})_{\phi}$ belong to the same Rouquier block of $(\mathcal{H}_{d,r})_\phi$.}
\end{apod}

The following result is a corollary of the above proposition. However, it can also be obtained independently using the Morita equivalences established by \cite{DiMa}:

\begin{proposition}\label{second characterization}
Let $\el$, $\mu$ be two $d$-partitions of $r$. The irreducible characters $(\chi_\el)_{\phi}$ and $(\chi_\mu)_{\phi}$ are in the same Rouquier block of $(\mathcal{H}_{d,r})_{\phi}$ if and only if  the following conditions are satisfied:
\begin{enumerate}
\item We have $\el^{(a)}=\mu^{(a)}$ for all $a \notin \{s,t\}$.
\item  If $\el^{st}:=(\el^{(s)},\el^{(t)})$, 
$\mu^{st}:=(\mu^{(s)},\mu^{(t)})$ and $l:=|\el^{st}|=|\mu^{st}|$, 
then the characters $(\chi_{\el^{st}})_{\vartheta}$ and $(\chi_{\mu^{st}})_{\vartheta}$ belong to the same Rouquier block of the cyclotomic Ariki-Koike algebra of $G(2,1,l)$ obtained via the specialization  
$$\vartheta: U_s\mapsto q^{m_s}, U_t \mapsto -q^{m_t}, X \mapsto q^n.$$
\end{enumerate}
\end{proposition}
\begin{apod}{Following \cite{DiMa}, Thm.1.1, we obtain that the algebra $(\mathcal{H}_{d,r})_\phi$ is Morita equivalent to the algebra
$$A:=\bigoplus_{\tiny
\begin{array}{c}
n_1,\ldots,n_{d-1} \geq 0\\
n_1+\ldots+n_{d-1}=r
\end{array}} (\mathcal{H}_{2,n_1})_{\phi'} \otimes \mathcal{H}(\mathfrak{S}_{n_2})_{\phi''}
\otimes\ldots\otimes  \mathcal{H}(\mathfrak{S}_{n_{d-1}})_{\phi''},
$$ where $\phi'$ is the restriction of $\phi$ to $\mathbb{Z}[u_s,u_s^{-1},u_t,u_t^{-1},x,x^{-1}]$
and $\phi''$ is the restriction of $\phi$ to $\mathbb{Z}[x,x^{-1}]$. Therefore, $(\mathcal{H}_{d,r})_\phi$ and $A$ have the same blocks.

Since $n \neq 0$, the Rouquier blocks of $\mathcal{H}(\mathfrak{S}_{n_2})_{\phi''}$,$\ldots$, 
$\mathcal{H}(\mathfrak{S}_{n_2})_{\phi''}$ are trivial. Thus we obtain that 
two irreducible characters $(\chi_\el)_{\phi}$ and $(\chi_\mu)_{\phi}$ are in the same Rouquier block of $(\mathcal{H}_{d,r})_{\phi}$ if and only if  the following conditions are satisfied:
\begin{enumerate}
\item We have $\el^{(a)}=\mu^{(a)}$ for all $a \notin \{s,t\}$.
\item If $\el^{st}:=(\el^{(s)},\el^{(t)})$, 
$\mu^{st}:=(\mu^{(s)},\mu^{(t)})$ and $l:=|\el^{st}|=|\mu^{st}|$, 
then the characters $(\chi_{\el^{st}})_{\phi'}$ and $(\chi_{\mu^{st}})_{\phi'}$ belong to the same block of
 $(\mathcal{H}_{2,l})_{\phi'}$ over the Rouquier ring of $\mathbb{Q}(\zeta_d)$.
\end{enumerate} 
Since the hyperplane  $kN+M_s-M_t=0$ is a $\mathfrak{p}$-essential hyperplane for $G(d,1,r)$, corollary $\ref{thm A}$ implies that the second condition is equivalent to saying that the $2$-partitions $\el^{st}$ and $\mu^{st}$ are
$\mathfrak{p}$-residue equivalent with respect to $\phi'$. By replacing  $\zeta_d^s q^{m_s}$ by $q^{m_s}$ and $\zeta_d^t q^{m_t}$ by $-q^{m_t}$ into the $\mathfrak{p}$-residues with respect to $\phi'$ of the nodes of $[\el^{st}]$ and $[\mu^{st}]$, we obtain their $2$-residues with respect to $\vartheta$. Therefore, the $2$-partitions $\el^{st}$ and $\mu^{st}$ are
$\mathfrak{p}$-residue equivalent with respect to $\phi'$ if and only if
they are $2$-residue equivalent with respect to $\vartheta$, \ie  the characters $(\chi_{\el^{st}})_{\vartheta}$ and 
$(\chi_{\mu^{st}})_{\vartheta}$ belong to the same Rouquier block of $(\mathcal{H}_{2,l})_{\vartheta}$.}
\end{apod}

\subsection*{\normalsize Case 3: Essential hyperplane $N=0$}

Let $$\phi : \left\{ 
\begin{array}{ll} 
u_j \mapsto \zeta_d^j q^{m_j}, (0 \leq j <d),\\ 
x \mapsto 1
\end{array} \right. 
$$
be a cyclotomic specialization associated with the essential hyperplane $$N=0.$$ Then we have the following result.

\begin{proposition}\label{essential hyperplane N=0}
Let $\lambda, \mu$ be two $d$-partitions of $r$.
The following assertions are equivalent:
\begin{description}
\item[(i)] The characters  $(\chi_\lambda)_\phi$ and  $(\chi_\mu)_\phi$ are in the same Rouquier block of $(\mathcal{H}_{d,r})_\phi$.
\item[(ii)]   $|\lambda^{(a)}|=|\mu^{(a)}|$ for all $a=0,1,\ldots,d-1$. 
\end{description}
\end{proposition}
\begin{apod}{(\textbf{(i) $\Rightarrow$ (ii)})  
Thanks to proposition $\ref{Rouquier blocks and central characters}$, we can assume that there exists a prime ideal $\mathfrak{p}$ of $\mathbb{Z}[\zeta_d]$ such that $(\chi_\lambda)_\phi$ and  $(\chi_\mu)_\phi$ belong to the same block of $k_\mathfrak{p}\mathcal{H}_\phi$ (where $k_\mathfrak{p}$ is the $\mathfrak{p}$-residue field of the Rouquier ring). Therefore, by theorem $\ref{thm residue}$, they must be $\mathfrak{p}$-residue equivalent with respect to $\phi$. Due to the form of the $\mathfrak{p}$-residue with respect to $\phi$ and the fact that the $m_a$ $(0 \leq a <d)$ can take any value, we must have 
$$\begin{array}{rcl}
|\lambda^{(a)}|=&\#\{(i,j,a)\,|\,  (1\leq i \leq h_\el^{(a)})(1 \leq j \leq \el_i^{(a)})\}&= \\
=& \#\{(i,j,a)\,|\,  (1\leq i \leq h_\mu^{(a)})(1 \leq j \leq \mu_i^{(a)})\}&=|\mu^{(a)}| 
\end{array}$$ 
for all $a=0,1,\ldots,d-1$.

(\textbf{(ii) $\Rightarrow$ (i)}) Let  $a \in \{0,1,\ldots,d-1\}$.
It is enough to show that if $\el $ and $\mu$ are two $d$-partitions of $r$ such that
\begin{center}
$|\el^{(a)}| =|\mu^{(a)}|$ and $\el^{(b)} = \mu^{(b)}$ for all $b \neq a$, 
\end{center}
then $(\chi_\el)_\phi$ and $(\chi_\mu)_\phi$ are in the same Rouquier block.

Set $l:=|\lambda^{(a)}|=|\mu^{(a)}|$. The generic Ariki-Koike algebra of the symmetric group 
$\mathfrak{S}_l$ specializes to the group algebra $\mathbb{Z}[\mathfrak{S}_l]$ when $x$ specializes to $1$ . 
For any finite group, it is well known that  $1$ is the only block-idempotent of the  group algebra over $\mathbb{Z}$ (see also \cite{Rou}, \S3, Rem.1). Thus, all irreducible characters of  $\mathfrak{S}_l$ belong to the same Rouquier block of $\mathbb{Z}[\mathfrak{S}_l]$. Corollary $\ref{thm A}$ implies that there exist a finite sequence of partitions of $l$, $\el_{(0)},\el_{(1)},\ldots,\el_{(m)}$ and a finite sequence of prime numbers of $\mathbb{Z}$, $p_1, p_2, \ldots, p_m$ such that
\begin{itemize}
\item $\el_{(0)}=\lambda^{(a)}$ and $\el_{(m)}=\mu^{(a)}$,
\item $\el_{(i-1)}$ and $\el_{(i)}$ are $(p_i)$-residue equivalent with respect to the specialization sending $x$ to $1$, for all $i=1,\ldots,m$.
\end{itemize}
We define $\el_{d,i}$ to be the $d$-partition of $r$ with
\begin{center}
$\el_{d,i}^{(a)}=\el_{(i)}$ and $\el_{d,i}^{(b)} = \el^{(b)}$ for all $b \neq a$.
\end{center} 
Let $\mathfrak{p}_i$ be a prime ideal of $\mathbb{Z}[\zeta_d]$ lying over the prime number $p_i$. Then we have
 \begin{itemize}
\item $\el_{d,0}=\lambda$ and $\el_{d,m}=\mu$,
\item $\el_{d,i-1}$ and $\el_{d,i}$ are $\mathfrak{p}_i$-residue equivalent with respect to $\phi$, for all $i=1,\ldots,m$.
\end{itemize}
Corollary $\ref{thm A}$ implies that $(\chi_\el)_\phi$ and $(\chi_\mu)_\phi$ are in the same Rouquier block of $(\mathcal{H}_{d,r})_\phi$.}
\end{apod}


\subsection*{\normalsize Conclusion}

Let $$\phi : \left\{ 
\begin{array}{ll} 
u_j \mapsto \zeta_d^j q^{m_j}, (0 \leq j <d),\\ 
x \mapsto q^n
\end{array} \right. 
$$
be a cyclotomic specialization for $\mathcal{H}_{d,r}$. Let $\el$ and $\mu$ be two $d$-partitions of $r$. We write $\el \sim_{R,\phi} \mu$ if there exist two integers $s$ and $t$ with $0 \leq s<t<d$ such that the following conditions are satisfied:
\begin{enumerate}
\item We have $\el^{(a)}=\mu^{(a)}$ for all $a \notin \{s,t\}$.
\item If $\el^{st}:=(\el^{(s)},\el^{(t)})$ and $\mu^{st}:=(\mu^{(s)},\mu^{(t)})$, then
$\mathrm{Contc}_{\el^{st}}= \mathrm{Contc}_{\mu^{st}}$ with respect to the weight system $(0,k)$, where $k$ is an integer such that $kn+m_s-m_t=0$
(or, equivalently, the irreducible characters $(\chi_{\el^{st}})_{\vartheta}$ and $(\chi_{\mu^{st}})_{\vartheta}$ belong to the same Rouquier block of the cyclotomic Ariki-Koike algebra of $G(2,1,l)$ obtained via the specialization  
$\vartheta: U_s\mapsto q^{m_s}, U_t \mapsto -q^{m_t}, X \mapsto q^n$).
\item There exists a prime ideal $\mathfrak{p}$ of $\mathbb{Z}[\zeta_d]$ such that 
$\zeta_d^s-\zeta_d^t \in \mathfrak{p}$.
\end{enumerate}
Thanks to propositions $\ref{no essential hyperplane}$, $\ref{essential hyperplane of type 1}$ and  $\ref{essential hyperplane N=0}$, proposition $\ref{explain AllBlocks}$ gives

\begin{theorem}\label{RB}
If $n \neq 0$, then two irreducible characters $(\chi_\el)_\phi$ and $(\chi_\mu)_\phi$ are in the same Rouquier block of 
$(\mathcal{H}_{d,r})_\phi$
if and only if there exists a finite sequence
$\el_{(0)},\el_{(1)},\ldots,\el_{(m)}$ of $d$-partitions of $r$  such that
\begin{itemize}
  \item $\el_{(0)}=\el$ and $\el_{(m)}=\mu$,
  \item for all $i$ $(1\leq i \leq m)$,\,\,
         we have $\el_{(i-1)} \sim_{R,\phi} \el_{(i)}$.
 \end{itemize}                   
 If $n = 0$, then two irreducible characters $(\chi_\el)_\phi$ and $(\chi_\mu)_\phi$ are in the same Rouquier block of 
$(\mathcal{H}_{d,r})_\phi$
if and only if there exists a finite sequence
$\el_{(0)},\el_{(1)},\ldots,\el_{(m)}$ of $d$-partitions of $r$  such that
\begin{itemize}
  \item $\el_{(0)}=\el$ and $\el_{(m)}=\mu$,
  \item for all $i$ $(1\leq i \leq m)$,\,\,
         we have $\el_{(i-1)} \sim_{R,\phi} \el_{(i)}$ or $|\el_{(i-1)}^{(a)}| = |\el_{(i)}^{(a)}|$ for all $a=0,1,\ldots,d-1$.
 \end{itemize}                   
\end{theorem} 

\subsection{The ``spetsial'' case}

In this section, we will show that the Rouquier blocks calculated by the algorithm of \cite{BK} are correct, when $\phi$ is the ``spetsial'' cyclotomic specialization (see example $\ref{spetsial}$). We are mostly interested in this case, because, as we have already mentioned, the Rouquier blocks of the ``spetsial'' cyclotomic Hecke algebra of a Weyl group coincide with its families of characters. \\

Let $$\phi : \left\{ 
\begin{array}{ll} 
u_j \mapsto \zeta_d^j q^{m_j}, (0 \leq j <d),\\ 
x \mapsto q
\end{array} \right. 
$$
be a cyclotomic specialization for $\mathcal{H}_{d,r}$. Let $\el$ and $\mu$ be two $d$-partitions of $r$. We write $\el \sim_{C,\phi} \mu$ if there exist two integers $s$ and $t$ with $0 \leq s<t<d$ such that the following conditions are satisfied:
\begin{enumerate}
\item We have $\el^{(a)}=\mu^{(a)}$ for all $a \notin \{s,t\}$.
\item If $\el^{st}:=(\el^{(s)},\el^{(t)})$ and $\mu^{st}:=(\mu^{(s)},\mu^{(t)})$, then
$\mathrm{Contc}_{\el^{st}}= \mathrm{Contc}_{\mu^{st}}$ with respect to the weight system $(m_s,m_t)$.
\end{enumerate}

\begin{proposition}\label{go to big prime}
Let $\el$ and $\mu$ be two $d$-partitions of $r$. We have that $\mathrm{Contc}_\el=
\mathrm{Contc}_\mu$ with respect to the weight system $(m_0,m_1,\ldots,m_{d-1})$ if and only if
there exists a finite sequence
$\el_{(0)},\el_{(1)},\ldots,\el_{(m)}$ of $d$-partitions of $r$  such that
\begin{itemize}
  \item $\el_{(0)}=\el$ and $\el_{(m)}=\mu$,
  \item for all $i$ $(1\leq i \leq m)$,\,\,
         we have $\el_{(i-1)} \sim_{C,\phi} \el_{(i)}$.
 \end{itemize}                   
\end{proposition}
\begin{apod}{We will first show that if $\el \sim_{C,\phi} \mu$, then $\mathrm{Contc}_\el=
\mathrm{Contc}_\mu$. Let $s,t$ be as in the definition of the relation $\sim_{C,\phi}$. 
Since $\mathrm{Contc}_{\el^{st}}= \mathrm{Contc}_{\mu^{st}}$ with respect to the weight system $(m_s,m_t)$, we have that $hc_{\el^{st}}=hc_{\mu^{st}}$. Moreover, $hc_{\el}^{(a)}=hc_\mu^{(a)}$ for all $a \neq s,t$. Therefore,
$$hc_\el=\mathrm{max}\{hc_{\el^{st}},(hc_{\el}^{(a)})_{a \neq s,t}\}=
                  \mathrm{max}\{hc_{\mu^{st}},(hc_{\mu}^{(a)})_{a \neq s,t}\}=hc_\mu.$$
Set $h:=hc_\el-hc_{\el^{st}}=hc_\mu-hc_{\mu^{st}}$. We have
$$Bc_\el^{(s)}=\beta_{\el}^{(s)}[hc_\el-h_{\el}^{(s)}+m_s]=\beta_{\el}^{(s)}[hc_{\el^{st}}-h_{\el}^{(s)}+m_s+h]=  Bc_{\el^{st}}^{(0)}[h].$$
Similarly, we obtain that
$$Bc_\el^{(t)}=Bc_{\el^{st}}^{(1)}[h],\,\, Bc_\mu^{(s)}=Bc_{\mu^{st}}^{(0)}[h] \,\textrm{ and }\, Bc_\mu^{(t)}=Bc_{\mu^{st}}^{(1)}[h].$$
Since
$$Bc_{\el^{st}}^{(0)} \cup Bc_{\el^{st}}^{(1)} =Bc_{\mu^{st}}^{(0)} \cup Bc_{\mu^{st}}^{(1)},$$
we have
$$Bc_{\el^{st}}^{(0)}[h] \cup Bc_{\el^{st}}^{(1)}[h] =Bc_{\mu^{st}}^{(0)}[h] \cup Bc_{\mu^{st}}^{(1)}[h]$$
and thus,
$$Bc_\el^{(s)} \cup Bc_\el^{(t)} =Bc_\mu^{(s)} \cup Bc_\mu^{(t)}.$$
Since $Bc_\el^{(a)}=Bc_\mu^{(a)}$ for all $a \neq s,t$, we deduce that $\mathrm{Contc}_\el=
\mathrm{Contc}_\mu$. 

Now let $\el$ and $\mu$ be two $d$-partitions of $r$ such that $\mathrm{Contc}_\el=\mathrm{Contc}_\mu$. 
Let $p$ be a prime number such that $p \geq d$.
We consider the cyclotomic specialization for $\mathcal{H}_{p,r}$
$$\bar{\phi} : \left\{ 
\begin{array}{ll} 
u_j \mapsto \zeta_p^j q^{m_j}, (0 \leq j <d),\\
u_i \mapsto \zeta_p^i q^{M}, (d \leq i <p),\\ 
x \mapsto q,
\end{array} \right.
$$
where $M > m_j+r$ for all $j\, ( 0 \leq j <d)$. We define the $p$-partition $\bar{\el}$ of $r$ as follows:
\begin{center}
$\bar{\el}^{(j)}:=\el^{(j)}$ for all $j\, ( 0 \leq j <d)$ and 
$\bar{\el}^{(i)}:=\emptyset$ for all $i\, ( d \leq i <p)$.
\end{center} 
Similarly, we define $\bar{\mu}$ as follows:
\begin{center}
$\bar{\mu}^{(j)}:=\mu^{(j)}$ for all $j\, ( 0 \leq j <d)$ and 
$\bar{\mu}^{(i)}:=\emptyset$ for all $i\, ( d \leq i <p)$.
\end{center} 

We have 
$hc_{\bar{\el}}^{(i)}=hc_{\bar{\mu}}^{(i)}=-M$  for all $i\, ( d \leq i <p)$. Moreover,
$hc_{\bar{\el}}^{(j)}>-M$ and $hc_{\bar{\mu}}^{(j)}>-M$  for all $j\, ( 0 \leq j <d)$.
Thus $hc_{\bar{\el}}=hc_\el=hc_\mu=hc_{\bar{\mu}}$. It is immediate, that
$\mathrm{Contc}_{\bar{\el}}=\mathrm{Contc}_{\bar{\mu}}$ with respect to the weight system
$(m_0,m_1,\ldots,m_{d-1},M,M,\ldots,M)$. 

 Since $p$ is a prime number, Theorem $\ref{BK}$ implies that the irreducible characters
$\chi_{\bar{\el}}$ and $\chi_{\bar{\mu}}$ belong to the same Rouquier block of $(\mathcal{H}_{p,r})_{\bar{\phi}}$. Due to Theorem $\ref{RB}$, there exists a finite sequence $\bar{\el}_{(0)},\bar{\el}_{(1)},\ldots,\bar{\el}_{(m)}$ of $p$-partitions of $r$  such that
\begin{itemize}
  \item $\bar{\el}_{(0)}=\bar{\el}$ and $\bar{\el}_{(m)}=\bar{\mu}$,
  \item for all $l$ $(1\leq l \leq m)$,\,\,
         we have $\bar{\el}_{(l-1)} \sim_{R,\bar{\phi}} \bar{\el}_{(l)}$ (and thus $\bar{\el}_{(l-1)} \sim_{C,\bar{\phi}} \bar{\el}_{(l)}$).
 \end{itemize} 
                  
Since $\bar{\el} \sim_{R,\bar{\phi}} \bar{\el}_{(1)}$ and $\bar{\phi}(x) =q \neq 1$, there exist two integers $s$ and $t$  with $0 \leq s<t<p$ such that $\bar{\el}$ and $\bar{\el}_{(1)}$ belong to the same Rouquier block associated with an essential hyperplane of the form 
$$kN+M_s-M_t=0, \textrm{ where } -r<k<r$$
and we have $k+m_s-m_t=0$. 
Since $M-m_j>r$ for all  $j\, ( 0 \leq j <d)$, we can't have $s < d \leq t$.
If $s \geq d$, then $\bar{\el}_{(1)}$ is a $p$-partition of $r$ if and only if $\bar{\el}_{(1)}=\bar{\el}$.
Thus, we must have $t<d$ and since ${\bar{\el}}^{(i)}=\emptyset$  for all $i\, ( d \leq i <p)$, we also have ${\bar{\el}_{(1)}}^{(i)}=\emptyset$ for all $i\, ( d \leq i <p)$.
Inductively, we obtain that ${\bar{\el}_{(l)}}^{(i)}=\emptyset$ for all $i\, ( d \leq i <p)$ and for all $l\, ( 1 \leq l \leq m)$ (the same result can be obtained from the fact that the charged content of two $p$-partitions linked by $\sim_{R,\bar{\phi}}$ is the same).
 
Let $l \in \{0,1,\ldots,m\}$. Let us define $\el_{(l)}$ to be the $d$-partition of $r$ such that
$\el_{(l)}^{(j)}:={\bar{\el}_{(l)}}^{(j)}$ for all $j\, ( 0 \leq j <d)$. Then we have
\begin{itemize}
  \item $\el_{(0)}=\el$ and $\el_{(m)}=\mu$,
  \item for all $l$ $(1\leq l \leq m)$,\,\,
         we have $\el_{(l-1)} \sim_{C,\phi} \el_{(l)}$.}
 \end{itemize}      
\end{apod}

Now assume that $\phi$ is the ``spetsial'' cyclotomic specialization, \ie
$$m_0=1 \,\textrm{ and }\, m_1=\ldots=m_{d-1}=0.$$
We are going to prove the following result.

\begin{proposition}\label{spetsial ok}
Let $\phi$ be the ``spetsial'' cyclotomic specialization. Two irreducible characters $(\chi_\el)_\phi$ and
$(\chi_\el)_\phi$ belong to the same Rouquier block of $(\mathcal{H}_{d,r})_\phi$ if and only if
$\mathrm{Contc}_\el=\mathrm{Contc}_\mu$. 
\end{proposition}
\begin{apod}{If $(\chi_\el)_\phi$ and
$(\chi_\el)_\phi$ belong to the same Rouquier block of $(\mathcal{H}_{d,r})_\phi$, then, by Theorem $\ref{BK}$, we have $\mathrm{Contc}_\el=\mathrm{Contc}_\mu$. 

Now let  $\el$ and $\mu$ be two $d$-partitions of $r$ such that $\mathrm{Contc}_\el=\mathrm{Contc}_\mu$. Thanks to proposition $\ref{go to big prime}$, we can assume that  $\el \sim_{C,\phi} \mu$. Then there exist  two integers $s$ and $t$ with $0 \leq s<t<d$ such that
$$\mathrm{Contc}_{\el^{st}}=\mathrm{Contc}_{\mu^{st}}    
\textrm{ and }\el^{(a)}=\mu^{(a)} \textrm{ for all } a\neq s,t.$$

Let us suppose that $d=p_1^{a_1}p_2^{a_2}\ldots p_n^{a_n}$, where $p_i$ are prime numbers such that $p_i \neq p_j$ for $i \neq j$. For $i=1,\ldots,n$, we set $c_i:=d/p_i^{a_i}$. 
Then $\mathrm{gcd}(c_i)=1$ and, by Bezout's theorem, there exist integers $(b_i)_{1 \leq i \leq n}$ such that
$\sum_{i=1}^n b_ic_i=1$. We have $s-t=\sum_{i=1}^n (s-t)b_ic_i$. We set $k_i:=(s-t)b_ic_i$ and we obtain that $s-t=\sum_{i=1}^n k_i$.

For all  $i=1,\ldots,n$, the element $1-\zeta_d^{c_i}$ belongs to the prime ideal of $\mathbb{Z}[\zeta_d]$ lying over the prime number $p_i$. So is $1-\zeta_d^{k_i}$.

Let $I$ be a subset of $\{1,\ldots,n\}$ minimal (with respect to inclusion) for the property
$$s-t \equiv \sum_{i \in I}k_i \,(\mathrm{mod}\,d),$$ 
\ie, if $J \subseteq I$ and
$$s-t \equiv \sum_{j \in J}k_j \,(\mathrm{mod}\,d),$$
then $J=I$. 
Without loss of generality, we can assume that $I=\{1,\ldots,n\}$. 
Now, for all $1 \leq m \leq n$, set
$$l_{m}:=\sum_{i=1}^mk_i \,(\mathrm{mod}\,d) \textrm{ and } l_0:=0.$$
Due to the minimality of $I$, we have $t+l_{i} \not\equiv s \,(\mathrm{mod}\,d)$ for all $i<n$.

The group $\mathfrak{S}_d$ acts naturally on the set of $d$-partitions of $r$: Let 
$\nu=(\nu^{(0)},\nu^{(1)},\ldots,\nu^{(d-1)})$ be a  $d$-partition of $r$.
If $\tau \in \mathfrak{S}_d$, then 
\begin{center}
$\tau(\nu)=(\nu^{(\tau(0))},\nu^{(\tau(1))},\ldots,\nu^{(\tau(d-1))})$.
 \end{center}
 For $a,b \in \{0,\ldots,d-1\}$, we denote by $(a,b)$ the corresponding transposition. If $a,b \neq 0$, then $\nu \sim_{C,\phi} (a,b)\nu$ (since the ordinary content is stable under the action of $(a,b)$).

For $i \in I$, set $\sigma_i:=(t+l_{i-1}  \,(\mathrm{mod}\,d),t+l_{i} \,(\mathrm{mod}\,d))$.
We have that the element
$$\zeta_d^{t+l_{i-1}}-\zeta_d^{t+l_{i}}=\zeta_d^{t+l_{i-1}}(1-\zeta_d^{k_i})$$
belongs to the prime ideal of $\mathbb{Z}[\zeta_d]$ lying over the prime number $p_i$. 
Therefore, if $t+l_{i-1} ,t+l_{i} \not\equiv 0 \,(\mathrm{mod}\,d)$ , then
$\nu \sim_{R,\phi} \sigma_i(\nu)$ for any $d$-partition $\nu$ of $r$.

Assume that $t+l_i \not\equiv 0 \,(\mathrm{mod}\,d)$ for all $i < n$. If $\sigma:=(t,t+l_{n-1} \,(\mathrm{mod}\,d))$, then
\begin{center}
$\sigma=\sigma_1 \circ \sigma_2 \circ \ldots \sigma_{n-2} \circ \sigma_{n-1} \circ \sigma_{n-2} \ldots
\circ \sigma_2 \circ \sigma_1.$
\end{center}
Theorem $\ref{RB}$ implies that $(\chi_\el)_\phi$ and $(\chi_{\sigma(\el)})_\phi$ belong to the same Rouquier block of $(\mathcal{H}_{d,r})_\phi$. The same holds for  $(\chi_\mu)_\phi$ and $(\chi_{\sigma(\mu)})_\phi$. Since $\el \sim_{C,\phi} \mu$ (with respect to $s$, $t$), we have that
$\sigma(\el) \sim_{C,\phi} \sigma(\mu)$ (with respect to $s$, $t+l_{n-1} \,(\mathrm{mod}\,d)$). Moreover, the element $\zeta_d^s-\zeta_d^{t+l_{n-1}}=\zeta_d^s(1-\zeta_d^{-k_n})$ belongs to the prime ideal of $\mathbb{Z}[\zeta_d]$ lying over the prime number $p_n$ and thus,  $\sigma(\el) \sim_{R,\phi} \sigma(\mu)$. Consequently, $(\chi_\el)_\phi$ and $(\chi_\mu)_\phi$ belong to the same Rouquier block of $(\mathcal{H}_{d,r})_\phi$. 

Now let us assume that there exists $1 \leq m < n$ such that
$$t+l_{i} \not\equiv0\,(\mathrm{mod}\,d)\textrm{ for all } i < m \textrm{ and } t+l_{m}\equiv 0\,(\mathrm{mod}\,d).$$ 
We will prove that $(\chi_\el)_\phi$ and $(\chi_\mu)_\phi$ belong to the same Rouquier block of $(\mathcal{H}_{d,r})_\phi$ by induction on $n-m$.

Let $m=n-1$. We have to distinguish two cases:
If  $k_{n-1} \not\equiv k_n \,(\mathrm{mod}\,d),$ then we have that $t+l_{n-2}+k_n \not\equiv0\,(\mathrm{mod}\,d)$ and we can rearrange the $k_i$ (exchanging $k_{n-1}$ and $k_n$) so that 
$t+l_{i} \not\equiv0\,(\mathrm{mod}\,d)\textrm{ for all } i < n$. This case has been covered above.

If  $k_{n-1} \equiv k_n \,(\mathrm{mod}\,d),$ we set
\begin{center}
$\sigma:=(t,t+l_{n-2}  \,(\mathrm{mod}\,d))=\sigma_1 \circ \sigma_2 \circ \ldots \sigma_{n-3} \circ \sigma_{n-2} \circ \sigma_{n-3} \ldots
\circ \sigma_2 \circ \sigma_1.$
\end{center}
Like above, we have that $(\chi_\el)_\phi$ and $(\chi_{\sigma(\el)})_\phi$ belong to the same Rouquier block of $(\mathcal{H}_{d,r})_\phi$. So do $(\chi_\mu)_\phi$ and $(\chi_{\sigma(\mu)})_\phi$.
Since the element
$$\zeta_d^s-\zeta_d^{t+l_{n-2}}=\zeta_d^s-\zeta_d^{s-k_{n-1}-k_n}=\zeta_d^s(1-\zeta_d^{-2k_n})$$
 belongs to the prime ideal of $\mathbb{Z}[\zeta_d]$ lying over the prime number $p_n$, we obtain that
 $\sigma(\el) \sim_{R,\phi} \sigma(\mu)$ and thus $(\chi_\el)_\phi$ and $(\chi_\mu)_\phi$ belong to the same Rouquier block of $(\mathcal{H}_{d,r})_\phi$.
 
Now assume that the result holds for integers greater than $m$. We will show that it holds for $m$:  Suppose that
$$t+l_{i} \not\equiv0\,(\mathrm{mod}\,d)\textrm{ for all } i < m \textrm{ and } t+l_{m}\equiv 0\,(\mathrm{mod}\,d).$$ 
 
 We again distinguish  two cases:
 If there exists $i_0 > m$ such that $k_{i_0} \not\equiv k_m \,(\mathrm{mod}\,d),$ then we have that $t+l_{m-1}+k_{i_0} \not\equiv 0\,(\mathrm{mod}\,d)$ and we can rearrange the $k_i$ (exchanging $k_{m}$ and $k_{i_0}$) so that $t+l_{i} \not\equiv0\,(\mathrm{mod}\,d)\textrm{ for all } i < m+1$. Now, the induction hypothesis and the case $t+l_{i} \not\equiv0\,(\mathrm{mod}\,d)\textrm{ for all } i < n$ cover all possibilities. Thus, the result is true.
 
 If $k_i \equiv k_m \,(\mathrm{mod}\,d),$ for all $i>m$, we set
 \begin{center}
$\sigma:=(t,t+l_{m-1}  \,(\mathrm{mod}\,d))=\sigma_1 \circ \sigma_2 \circ \ldots \sigma_{m-2} \circ \sigma_{m-1} \circ \sigma_{m-2} \ldots
\circ \sigma_2 \circ \sigma_1.$
\end{center}
Again we have that $(\chi_\el)_\phi$ and $(\chi_{\sigma(\el)})_\phi$ belong to the same Rouquier block of $(\mathcal{H}_{d,r})_\phi$. So do $(\chi_\mu)_\phi$ and $(\chi_{\sigma(\mu)})_\phi$.
Since the element
$$\zeta_d^s-\zeta_d^{t+l_{m-1}}=\zeta_d^{t+l_n}-\zeta_d^{t+l_{m-1}}=
\zeta_d^{t+l_{m-1}}(\zeta_d^{(n-m+1)k_m}-1)$$
 belongs to the prime ideal of $\mathbb{Z}[\zeta_d]$ lying over the prime number $p_m$, we obtain that
 $\sigma(\el) \sim_{R,\phi} \sigma(\mu)$ and thus $(\chi_\el)_\phi$ and $(\chi_\mu)_\phi$ belong to the same Rouquier block of $(\mathcal{H}_{d,r})_\phi$.}
\end{apod}

\subsection{Functions $a$ and $A$}

Let $$\phi : \left\{ 
\begin{array}{ll} 
u_j \mapsto \zeta_d^j q^{m_j}, (0 \leq j <d),\\ 
x \mapsto q^n
\end{array} \right. 
$$
be a cyclotomic specialization for $\mathcal{H}_{d,r}$. If $n \neq 0$, then \cite{BK}, Prop.3.18 implies that the functions $a$ and $A$ (see section $2.4$) are constant on the Rouquier blocks of $(\mathcal{H}_{d,r})_\phi$. We will show that this is also true for $n=0$. Thanks to theorem $\ref{RB}$, it suffices to show that

\begin{proposition}\label{a}
Let $\el$ and $\mu$ be two $d$-partitions of $r$. Let $\phi$ be a cyclotomic specialization associated with the essential hyperplane $N=0$. If $(\chi_\el)_\phi$ and $(\chi_\mu)_\phi$ belong to the same Rouquier block of $(\mathcal{H}_{d,r})_\phi$, then 
$$a((\chi_\el)_\phi)=a((\chi_\mu)_\phi) \textrm{ and } A((\chi_\el)_\phi)=A((\chi_\mu)_\phi).$$
\end{proposition}
\begin{apod}{Thanks to proposition $\ref{aA}$, we have that 
$$a((\chi_\el)_\phi)+A((\chi_\el)_\phi) = a((\chi_\mu)_\phi)+A((\chi_\mu)_\phi).$$
Thus, it is enough to show that $A((\chi_\el)_\phi)=A((\chi_\mu)_\phi).$

Set $L:=\mathrm{max}\{h_\el,h_\mu\}$. Using the notations of proposition $\ref{schur}$, it is straightforward to check that, for $x=1$, the term $\delta_\el$ doesn't depend on the $d$-partition $\el$. Consequently, we obtain that
$A((\chi_\el)_\phi)=A((\chi_\mu)_\phi)$ if and only if
$$
\begin{array}{l}
\mathrm{deg}_q(\prod_{0 \leq s,t <d}\prod_{b_s \in B^{(s)}_{\el,L}}\prod_{1 \leq k \leq b_s} 
(\zeta_d^sq^{m_s}-\zeta_d^tq^{m_t}))=\\ \mathrm{deg}_q(\prod_{0 \leq s,t <d}\prod_{b_s \in B^{(s)}_{\mu,L}}\prod_{1 \leq k \leq b_s} 
(\zeta_d^sq^{m_s}-\zeta_d^tq^{m_t})).
\end{array}$$
Set
\begin{center}
$f_\el(q):=\prod_{0 \leq s,t <d}\prod_{b_s \in B^{(s)}_{\el,L}}\prod_{1 \leq k \leq b_s} 
(\zeta_d^sq^{m_s}-\zeta_d^tq^{m_t}).$
\end{center}
We have that
$$
\begin{array}{rcl}
f_\el(q)
&=&\prod_{0 \leq s,t <d}\prod_{b_s \in B^{(s)}_{\el,L}}((\zeta_d^sq^{m_s}-\zeta_d^tq^{m_t})^{b_s})\\
&=&\prod_{0 \leq s,t <d}((\zeta_d^sq^{m_s}-\zeta_d^tq^{m_t})^{\sum b_s})\\
&=&\prod_{0 \leq s,t <d}((\zeta_d^sq^{m_s}-\zeta_d^tq^{m_t})^{|\el^{(s)}|+\binom{L}{2}}))
\end{array}$$

Since $(\chi_\el)_\phi$ and $(\chi_\mu)_\phi$ belong to the same Rouquier block of $(\mathcal{H}_{d,r})_\phi$, by theorem $\ref{essential hyperplane N=0}$, we have $|\lambda^{(s)}|=|\mu^{(s)}|$ for all $s=0,1,\ldots,d-1$. Thus, $f_\el(q)=f_\mu(q)$, which implies that
$\mathrm{deg}_q(f_\el(q))=\mathrm{deg}_q(f_\mu(q))$. Therefore, we obtain that $A((\chi_\el)_\phi)=A((\chi_\mu)_\phi)$.}
\end{apod}


\begin{thebibliography}{00}
\addcontentsline{toc}{chapter}{Bibliography}
\bibitem{Ar1} S. Ariki, \emph{On the semi-simplicity of the Hecke algebra of $(\mathbb{Z}/r\mathbb{Z}) \wr \mathfrak{S}_n$}, J. Algebra 169, No. 1 (1994), 216-225. 
\bibitem{ArKo} S. Ariki, K. Koike, \emph{A Hecke algebra of $(\mathbb{Z}/r\mathbb{Z}) \wr S_n$ and construction of its irreducible representations}, Adv. in Math. 106 (1994), 216-243. 
\bibitem{Bro}
M. Brou{\'e}, \emph{On representations of symmetric algebras: An
introduction}, Notes by Markus Stricker, Forschungsinstitut f{\"u}r
Mathematik, ETH, Z{\"u}rich, 1991.
\bibitem{BK}  M. Brou{\'e}, S. Kim, \emph{Familles de caract{\`e}res des alg{\`e}bres de Hecke
cyclotomiques}, Adv. in Mathematics 172 (2002), 53-136.
\bibitem{BM} M. Brou{\'e}, G. Malle,  \emph{Zyklotomische Heckealgebren}, Ast\'{e}risque 212 (1993), 119-189.
\bibitem{BMM2}  M. Brou{\'e}, G. Malle, J. Michel, \emph{Towards Spetses I}, Trans. Groups 4, No. 2-3 (1999), 157-218.
\bibitem{BMR} M. Brou{\'e}, G. Malle, R. Rouquier, \emph{Complex
reflection groups, braid groups, Hecke algebras}, J. reine angew.
Math. 500 (1998), 127-190.
\bibitem{ChDeg} M. Chlouveraki, \emph{Degree and valuation of the Schur elements of cyclotomic Hecke algebras}, arXiv:0802.4306.
\bibitem{Chlou} M. Chlouveraki, \emph{On the cyclotomic Hecke algebras of complex reflection groups},
Ph.D. thesis, Universit{\'e} Paris 7, 2007 (available online at arXiv:0710.0776v1). 
\bibitem{DiMa} R. Dipper, A. Mathas, \emph{Morita Equivalences of Ariki-Koike algebras}, Math. Z. 240 (2002), 579-610.
\bibitem{Ge} M. Geck, \emph{Beitr{\"a}ge zur Darstellungstheorie von Iwahori-Hecke-Algebren}, RWTH Aachen, Habilitations-schrift, 1993.
\bibitem{GePf} M. Geck, G. Pfeiffer, Characters of Coxeter groups and Iwahori-Hecke algebras,
LMS monographs, New series no. 21, Oxford Univ. Press, 2000.
\bibitem{GeRo} M. Geck, R. Rouquier, \emph{Centers and simple modules for Iwahori-Hecke algebras}, Progress in Math. 141, Birkha{\"u}ser (1997), 251-272.
\bibitem{GIM}  M. Geck, L. Iancu, G. Malle, \emph{Weights of Markov traces and generic degrees}, Indag. Math. 11 (2000), 379-397.
\bibitem{GrLe} J. Graham, G. Lehrer, \emph{Cellular algebras}, Invent. Math. 123 (1996), 1-34. 
\bibitem{Gy} A. Gyoja, \emph{Cells and modular representations of
Hecke algebras}, Osaka J. Math. 33 (1996), 307-341.
\bibitem{Kim} S. Kim, \emph{Families of the characters of the cyclotomic Hecke algebras of}
$G(de,e,r)$, J. Algebra 289 (2005), 346-364.
\bibitem{Lu1} G. Lusztig, Characters of Reductive Groups
over a Finite Field, Annals of Mathematical Studies, Vol. 107,
Princeton Univ. Press, Princeton, NJ, 1984.
\bibitem{Lu2} G. Lusztig,  \emph{Leading coefficients of character values of Hecke algebras}, Proc. Symp. Pure Math., vol. 47(2), Amer. Math. Soc., Providence, RI, 1987, 235-262.
\bibitem{LyMa} S. Lyle, A. Mathas, \emph{Blocks of cyclotomic Hecke algebras}, Adv. Math., 216 (2007), 854-878. 
\bibitem{Ma4}  G. Malle, \emph{On the rationality and fake degrees of characters of cyclotomic algebras}, J. Math. Sci. Univ.
Tokyo 6 (1999), 647-677.
\bibitem{MaRo} G. Malle, R. Rouquier, \emph{Familles de caract{\`e}res de groupes de r{\'e}flexions complexes}, Representation theory 7 (2003), 610-640.
\bibitem{Mat}  A. Mathas,  \emph{Matrix units and generic degrees for the Ariki-Koike algebras}, J. Algebra 281 (2004), 695-730.
\bibitem{Na} M. Nagata, Local rings, Interscience tracts in pure and applied mathematics no.13, Interscience
publishers, U.S.A., 1962.
\bibitem{Rou}  R. Rouquier, \emph{Familles et blocs d'alg{\`e}bres de Hecke}, C. R. Acad. Sciences 329 (1999), 1037-1042.
\bibitem{The}  J. Th{\'e}venaz, $G$-algebras and Modular Representation Theory, Oxford Science Publications, Clarendon Press, 1995.







\end{thebibliography}
\end{document}